\documentclass[12pt,leqno]{article}
\usepackage{amsmath, amssymb, indentfirst}

\newcommand{\sect}[1]{\setcounter{equation}{0}\section{#1}}

\def\epsilon{\varepsilon}

\begin{document}

\LARGE 
\large

\noindent 
{\bf Heat equation and Schr\"{o}dinger equation}

\noindent 
{\bf with translation invariance}

\noindent  
{\bf on the infinite-dimensional vector space $\mathbb R^\infty$}

\large
\normalsize

\vspace*{0.8em}

\noindent

\hfill Hiroki Yagisita (Kyoto Sangyo University) 

\vspace*{3.2em} 

\normalsize

The standard Laplacian $-\triangle_{\mathbb R^n}$ in $L^2(\mathbb R^n)$ 
is self-adjoint and translation invariant on the finite-dimensional linear space $\mathbb R^n$. 
In this paper, we define a translation invariant operator 
$-\triangle_{\mathbb R^\infty}$ on $\mathbb R^\infty$ 
as a non-negative self-adjoint operator in some non-separable Hilbert space $L^2(\mathbb R^\infty)$. 
The set $L^2(\mathbb R^\infty)$ is a translation invariant subset of the set $CM(\mathbb R^\infty)$ 
of all complex measures on the product measurable space $\mathbb R^\infty$.  
Furthermore, we show that for any $f\in L^2(\mathbb R^n)$ 
and any $u\in L^2(\mathbb R^\infty)$, the separations of variables 
$e^{\triangle_{\mathbb R^\infty}t}(f\otimes u)
=(e^{\triangle_{\mathbb R^n}t}f)\otimes 
(e^{\triangle_{\mathbb R^\infty}t}u) \ (t\in [0,+\infty))$ 
and $e^{\sqrt{-1}\triangle_{\mathbb R^\infty}t}(f\otimes u)
=(e^{\sqrt{-1}\triangle_{\mathbb R^n}t}f)\otimes 
(e^{\sqrt{-1}\triangle_{\mathbb R^\infty}t}u) \ (t\in (-\infty,+\infty))$
hold. This clearly shows that $-\triangle_{\mathbb R^\infty}$ 
is an analog of $-\triangle_{\mathbb R^n}$. 

The starting point for the discussion in this paper 
is to naturally introduce a translation invariant structure of Hilbert space into $CM(\mathbb R^\infty)$. 
$L^2(\mathbb R^\infty)$ is a closed linear subspace of $CM(\mathbb R^\infty)$.  
The inner product of $L^2(\mathbb R^\infty)$ is defined as that of $CM(\mathbb R^\infty)$.
For a manifold, H\"{o}rmander defined an inner product that does not depend on a particular measure. 
In fact, the way we introduce the inner product into $CM(\mathbb R^\infty)$ is a generalization of his. 
Not only is a statistical manifold on $\mathbb R^\infty$ a submanifold of $CM(\mathbb R^\infty)$, 
but the real inner product $\mathrm{Re}(\langle \cdot, \cdot \rangle_{CM(\mathbb R^\infty)})$ 
induces Fisher information metric.

\vfill









\newpage

\sect{Introduction}
In this paper, we define a translation invariant operator 
$-\triangle_{\mathbb R^\infty}$ on the infinite-dimensional vector space $\mathbb R^\infty$ 
as a non-negative self-adjoint operator and we examine evolution equations 
$u_t=\triangle_{\mathbb R^\infty}u$ and $u_t=\sqrt{-1}\triangle_{\mathbb R^\infty}u$ 
to show that $-\triangle_{\mathbb R^\infty}$ 
is an analog of the standard Laplacian $-\triangle_{\mathbb R^n}$ in $L^2(\mathbb R^n)$.  
Our construction of $-\triangle_{\mathbb R^\infty}$ 
is achieved by combining well basic results 
widely used in {\bf finite-dimensional} analysis (e.g., [2, 3, 5, 7, 8, 12]), 
although it is an infinite-dimensional object. 
On the other hand, our method 
cannot be immediately applied to constructing an analog on a smooth domain of $\mathbb R^\infty$, 
so we look forward to further research in the future. 
There is no known previous research that has given an analog on an infinite-dimensional linear space 
that is self-adjoint {\bf and} translation invariant.  
In this sense, there does not seem to be any clear related literature. 
On the other hand, the amount of previous research 
on things that seem in some sense to be analogs on infinite-dimensional linear spaces 
that are self-adjoint {\bf or} translation invariant seems to be vast. 
We find it difficult to provide an unbiased citation. 

Let the set of all measurable sets of $\mathbb R$ 
be the topological $\sigma$-algebra of $\mathbb R$ 
(i.e., the smallest $\sigma$-algebra on $\mathbb R$ containing all open sets of $\mathbb R$).  
Let $\mathbb R^\infty$ denote the countable product measurable space $\prod_{n\in\mathbb N}\mathbb R$.  
In Sections 2, 3 and 4, we define some Hilbert space $L^2(\mathbb R^\infty)$ 
which is a subset of the set of all complex measures on $\mathbb R^\infty$  
and a non-negative self-adjoint operator $-\triangle_{\mathbb R^\infty}$ in $L^2(\mathbb R^\infty)$. 
In Section 5, we set the stage for the following sections. 
In Section 6, we examine evolution equations 
$u_t=\triangle_{\mathbb R^\infty}u$ and $u_t=\sqrt{-1}\triangle_{\mathbb R^\infty}u$ 
to show that $-\triangle_{\mathbb R^\infty}$ 
is an analog of the standard Laplacian $-\triangle_{\mathbb R^n}$ in $L^2(\mathbb R^n)$. 
In Section 7, we show that $L^2(\mathbb R^\infty)$ is not separable.  
In Section 8, we show that $L^2(\mathbb R^\infty)$ and $\triangle_{\mathbb R^\infty}$ 
are translation invariant on $\mathbb R^\infty$. 

More specifically, it is as follows.

There does not exist a $\sigma$-finite measure $\mu$ on $\mathbb R^\infty$ with translation invariance 
that satisfies $\mu([0,1)^\infty)=1$ (i.e., ideal Lebesgue measure on $\mathbb R^\infty$). 
Therefore, the space of square-integrable functions on $\mathbb R^\infty$ seems indefinable. 
However, on the other hand, according to Born and Heisenberg probabilistic interpretation 
of quantum mechanical wavefunction, when $\Omega$ is a measurable space, 
for a measure $\mu$ on $\Omega$ and a function $f\in L^2(\mu)$ 
that satisfy $\|f\|_{L^2(\mu)}\not=0$, the probabilistic interpretation of the state vector $f$ 
for position measurement is the probability measure 
$$\frac{1}{\|f\|_{L^2(\mu)}^2}|f|^2d\mu$$
on $\Omega$. This probability measure 
is the normalization of the total variation of the complex measure 
$$f|f|d\mu$$ 
on $\Omega$. So, it raises the following question. 
For a measurable space $\Omega$, let $CM(\Omega)$ denote the set of all complex measures on $\Omega$, 
into $CM(\Omega)$, is it possible to introduce the structure of Hilbert space 
such that for any measure $\mu$ on $\Omega$ 
and any $f_1, f_2\in L^2(\mu)$, 
$$\langle f_1|f_1|d\mu, f_2|f_2|d\mu \rangle_{CM(\Omega)}=\int_\Omega f_1\overline{f_2}d\mu$$ 
holds ? As a matter of fact, 
given restricting $\mu$ to be $\sigma$-finite, 
this can be easily done consistently using Radon-Nikodym theorem, 
as we will do in Section 2. 
That is, for any measurable space $\Omega$, 
the structure of Hilbert space such that for any $\sigma$-finite measure $\mu$ on $\Omega$ 
and any $f_1, f_2\in L^2(\mu)$, 
$\langle f_1|f_1|d\mu, f_2|f_2|d\mu \rangle_{CM(\Omega)}=\int_\Omega f_1\overline{f_2}d\mu$ 
holds is introduced into $CM(\Omega)$. 
In particular, when $\Omega$ is a measurable space, 
for any $\sigma$-finite measure $\mu$ on $\Omega$, 
the isometry 
$$f \ \ \ \mapsto \ \ \ f|f|d\mu$$
from $L^2(\mu)$ into $CM(\Omega)$ is determined. 
So,  we can think that the well-known Hilbert space $L^2(\mathbb R^n)$ 
is a closed linear subspace of Hilbert space $CM(\mathbb R^n)$.  
Since the sum in $CM(\Omega)$ that fits the inner product is defined in Definition 4 and it is simple, 
readers who have doubts here should take a look at Definition 4.

For a vector $a\in \mathbb R^\infty$, let $\tau_a$ denote the translation on $CM(\mathbb R^\infty)$ by the vector $a$. 
Let $e_1:=(1,0,0,\cdots) \in \mathbb R^\infty, e_2:=(0,1,0,\cdots) \in \mathbb R^\infty, \cdots$. 
Let $L^2_k(\mathbb R^\infty)$ denote the set of all $u\in CM(\mathbb R^\infty)$ such that 
$\lim_{h\downarrow+0}\|\tau_{he_k}u-u\|_{CM(\mathbb R^\infty)}=0$ holds. 
Let $H^1_k(\mathbb R^\infty)$ denote the set of all $u\in CM(\mathbb R^\infty)$ 
such that there uniquely exists $v\in CM(\mathbb R^\infty)$ such that 
$\lim_{h\downarrow+0}\|\frac{\tau_{he_k}u-u}{h}-v\|_{CM(\mathbb R^\infty)}=0$ holds. 
For $u\in H^1_k(\mathbb R^\infty)$, let $\frac{\partial u}{\partial x_k}$ denote $-v$. In Section 3, we show 
that $\tau_a$ is a unitary operator in $CM(\mathbb R^\infty)$. 
So, by Stone theorem, $\sqrt{-1}\frac{\partial}{\partial x_k}$ is a self-adjoint operator in $L^2_k(\mathbb R^\infty)$.

In Section 4, from the very simple quadratic form (Hermitian form) 
$$\sum_{k\in \mathbb N}\langle \sqrt{-1}\frac{\partial u_1}{\partial x_k}, 
\sqrt{-1}\frac{\partial u_2}{\partial x_k}\rangle_{CM(\mathbb R^\infty)},$$ 
we define a closed linear subspace $L^2(\mathbb R^\infty)$ of $CM(\mathbb R^\infty)$ 
and a non-negative self-adjoint operator $-\triangle_{\mathbb R^\infty}$ in $L^2(\mathbb R^\infty)$. 
First, we introduce Sobolev type space $H^1(\mathbb R^\infty)$ in $CM(\mathbb R^\infty)$ 
by the inner product 
$$\langle u_1, u_2 \rangle_{H^1(\mathbb R^\infty)}
:=\langle u_1, u_2 \rangle_{CM(\mathbb R^\infty)}
+\sum_{k\in\mathbb N}
\langle\sqrt{-1}\frac{\partial u_1}{\partial x_k},
\sqrt{-1}\frac{\partial u_2}{\partial x_k}\rangle_{CM(\mathbb R^\infty)}.$$
The closure of $H^1(\mathbb R^\infty)$ in $CM(\mathbb R^\infty)$ is denoted by $L^2(\mathbb R^\infty)$. 
We show that $H^1(\mathbb R^\infty)$ is Hilbert space. 
Therefore, by $|\langle f, \varphi \rangle_{CM(\mathbb R^\infty)}|\leq \|f\|_{CM(\mathbb R^\infty)}\|\varphi\|_{H^1(\mathbb R^\infty)}$ 
and Riesz theorem, for any $f\in L^2(\mathbb R^\infty)$, there uniquely exists $u\in H^1(\mathbb R^\infty)$ 
such that for any $\varphi \in H^1(\mathbb R^\infty)$, 
$$\langle f, \varphi \rangle_{CM(\mathbb R^\infty)}=\langle u, \varphi \rangle_{H^1(\mathbb R^\infty)}$$
holds. This map $f \mapsto u$ from $L^2(\mathbb R^\infty)$ to $H^1(\mathbb R^\infty)$ 
is denoted by $(1-\triangle_{\mathbb R^\infty})^{-1}$. 
We show that $(1-\triangle_{\mathbb R^\infty})^{-1}$ is self-ajoint in $L^2(\mathbb R^\infty)$ and injective.  
So, the operator $\triangle_{\mathbb R^\infty}$ in $L^2(\mathbb R^\infty)$ defined by  
$$\triangle_{\mathbb R^\infty}u:=u-((1-\triangle_{\mathbb R^\infty})^{-1})^{-1}u$$
is self-ajoint.  We show that $-\triangle_{\mathbb R^\infty}$ is non-negative.

In Section 5, when $\Omega_1$ and $\Omega_2$ are measurable spaces, 
$u_1$ is a complex measure on $\Omega_1$ 
and $u_2$ is a complex measure on $\Omega_2$, it 
is shown that the product $u_1\cdot u_2$ 
can be naturally defined as a complex measure on $\Omega_1\times\Omega_2$ 
and $\|u_1\cdot u_2\|_{CM(\Omega_1\times\Omega_2)}=\|u_1\|_{CM(\Omega_1)}\|u_2\|_{CM(\Omega_2)}$ holds. 
Since the product is defined in Definition 26 and it is simple, 
readers who have doubts here should take a look at Definition 26. 

For $f\in L^2(\mathbb R^n)$ and $u\in CM(\mathbb R^\infty)$, 
let $f\otimes u\in CM(\mathbb R^\infty)$ be defined as 
$$(f\otimes u)(x_1,x_2,\cdots)$$
$$:=(f(x_1,x_2,\cdots,x_n)|f(x_1,x_2,\cdots,x_n)|dx_1dx_2\cdots dx_n)
\cdot u(x_{n+1},x_{n+2},\cdots).$$ 
In Section 6, we show that for any $f_0\in L^2(\mathbb R^n)$ 
and any $u_0\in L^2(\mathbb R^\infty)$, the formulas for separation of variables 
$$e^{\triangle_{\mathbb R^\infty}t}(f_0\otimes u_0)
=(e^{\triangle_{\mathbb R^n}t}f_0)\otimes 
(e^{\triangle_{\mathbb R^\infty}t}u_0) \ \ \ (t\in [0,+\infty)),$$ 
$$e^{\sqrt{-1}\triangle_{\mathbb R^\infty}t}(f_0\otimes u_0)
=(e^{\sqrt{-1}\triangle_{\mathbb R^n}t}f_0)\otimes 
(e^{\sqrt{-1}\triangle_{\mathbb R^\infty}t}u_0) \ \ \ (t\in (-\infty,+\infty))$$
hold. This clearly shows that $\triangle_{\mathbb R^\infty}$ 
is an analog of the standard Laplacian $\triangle_{\mathbb R^n}$. 
In order to prove this, we had to resort to technical ingenuity. 
The proof path for the separations of variables is not very clear even to the author himself, 
so we hope that readers will improve it. 

In Section 7, we show that $L^2(\mathbb R^\infty)$ has an uncountable orthogonal system
 (i.e., $L^2(\mathbb R^\infty)$ is not separable). 

In Section 8, just to make sure, we 
confirm that $L^2(\mathbb R^\infty)$ and $\triangle_{\mathbb R^\infty}$ 
are translation invariant on $\mathbb R^\infty$.

Section 6, Section 7 and Section 8 can each be read independently. 
In addition, as appropriate, we once again state the proof of some fundamental facts that should hold, 
because there does not seem to be much accessible well-known literature 
that explicitly states the proof (and the author is the kind of person who worries about such things). 
Readers who are not too concerned can feel free to ignore such statements. 
On the one hand, in a familiar way, it uses well-known basic theorems such as Radon-Nikodym theorem 
without any specific reference to them.

\newpage 

\sect{Square root of density}
Throughout this paper, we may use the following three simple facts 
(Lemma 1, Lemma 2 and Remark) without specific mention. 
In a familiar way, we also use Radon-Nikodym theorem without any specific reference. 

{\bf Lemma 1:} 
Let $\eta$ be the map from $\mathbb C$ to $\mathbb C$ defined for $z\in \mathbb C$, 
as $\eta(z):=z|z|$. 
Let $\zeta$ be the map from $\mathbb C$ to $\mathbb C$ defined for $w\in \mathbb C$, 
as $\zeta(w):=w|w|^{-\frac{1}{2}}$ when $w\not=0$ holds 
and $\zeta(w):=0$ when $w=0$ holds.  Then, $\eta$ and $\zeta$ are continuous. 
$\zeta\circ\eta$ and $\eta\circ\zeta$ are the identity map. 
\hfill 
$\square$

{\bf Lemma 2:} 
Let $\Omega$ be a measurable space. 
Let $\{\mu_n\}_{n=1}^\infty$ be a sequence of $\sigma$-finite measures on $\Omega$.  
Let $\{u_n\}_{n=1}^\infty$ be a sequence of complex measures on $\Omega$. 
Then, there exists a finite measure $\nu$ on $\Omega$ such that 
for any $n\in \mathbb N$, $\mu_n$ and $u_n$ are absolutely continuous 
with respect to $\nu$. 

{\bf Proof:}  
There exists $\{E_{n,m}\}_{n,m=1}^\infty$ such that for any $n$, 
$\Omega=\sqcup_{m=1}^\infty E_{n,m}$ holds and for any $m$, 
$\mu_n(E_{n,m})<+\infty$ holds.  
As $|u_n|$ is the total variation of $u_n$, 
let
$$\nu(E):=\sum_n\left(\frac{1}{2^n}\left(\frac{|u_n|(E)}{1+|u_n|(\Omega)}
+\sum_m\left(\frac{1}{2^m}\frac{\mu_n(E\cap E_{n,m})}{1+\mu_n(E_{n,m})}\right)\right)\right).$$ 
\hfill 
$\blacksquare$

{\bf Remark:} 
Let $\Omega$ be a measurable space. 
Let $\mu$ be a measure on $\Omega$. 
Let $\rho$ be a $[0,+\infty)$-valued measurable function on $\Omega$. 
Let $f$ be a complex measurable function on $\Omega$. 
Then, the followings hold. 

(1) Suppose that $f$ is non-negative. 
Then, $f(\rho d\mu)=(f\rho)d\mu$ holds. 

(2) Suppose that $f\in L^1(\rho d\mu)$ or $f\rho\in L^1(\mu)$ holds. 
Then, $f\in L^1(\rho d\mu)$, $f\rho\in L^1(\mu)$ and $f(\rho d\mu)=(f\rho)d\mu$ hold. 

{\bf Proof:} 
Although it is a natural result, we include the proof just in case. 

(1) There exists a monotonic non-decreasing sequence $\{g_n\}_n$ 
of non-negative simple measurable functions 
such that for any $x$, $\lim_n |g_n(x)-f(x)|=0$ holds. 
Then, for any measurable set $E$, 
$(f(\rho d\mu))(E)
=\int_E f(\rho d\mu)
=\lim_n (\int_E g_n(\rho d\mu))
=\lim_n (\int_E (g_n\rho)d\mu)
=\int_E (f\rho)d\mu
=((f\rho)d\mu)(E)$ 
holds. 

(2) From (1), $f\in L^1(\rho d\mu)$ and $f\rho\in L^1(\mu)$ hold. 
So, there exist $h_1, h_2, h_3, h_4\in L^1(\rho d\mu)$ such that $h_1, h_2, h_3$ and $h_4$ 
are non-negative and 
$$f=(h_1-h_2)+\sqrt{-1}(h_3-h_4)$$ 
holds. Then, from (1), for any measurable set $E$, 
$(f(\rho d\mu))(E)
=\int_E f(\rho d\mu)
=(\int_Eh_1(\rho d\mu)-\int_Eh_2(\rho d\mu))
+\sqrt{-1}(\int_Eh_3(\rho d\mu)-\int_Eh_4(\rho d\mu))
=(\int_E(h_1\rho)d\mu-\int_E(h_2\rho)d\mu)
+\sqrt{-1}(\int_E(h_3\rho)d\mu-\int_E(h_4\rho)d\mu)
=\int_E(f\rho)d\mu
=((f\rho)d\mu)(E)$ holds. 
\hfill 
$\blacksquare$

{\bf Definition 3 ($CM(\Omega)$):} 
Let $\Omega$ be a measurable space. Then, let $CM(\Omega)$ denote the set 
of all complex measures on $\Omega$. 
\hfill 
$\square$

{\bf Definition 4 (sum):} 
Let $\Omega$ be a measurable space. Let $u_1, u_2 \in CM(\Omega)$.
Then, there uniquely exists $v\in CM(\Omega)$ such that 
the following holds.  
There exist a $\sigma$-finite measure $\mu$ on $\Omega$ 
and $f_1, f_2\in L^2(\mu)$ such that 
$u_1=f_1|f_1|d\mu$, $u_2=f_2|f_2|d\mu$ 
and $v=(f_1+f_2)|f_1+f_2|d\mu$ hold.
Let $u_1+u_2 \in CM(\Omega)$ be defined as $u_1+u_2:=v$. 

{\bf Proof:} 
From Lemma 1 and Lemma 2, existence is easy. We show uniqueness. 
Suppose that $(v_1,\mu_1,f_{1,1},f_{2,1})$ and 
$(v_2,\mu_2,f_{1,2},f_{2,2})$ each satisfy the condition of the definition. 
Then, from Lemma 2, there exist a finite measure $\nu$ 
and $[0,+\infty)$-valued measurable functions $\rho_1, \rho_2$ such that 
$\mu_1=\rho_1d\nu$ and $\mu_2=\rho_2d\nu$ hold. So, 
$u_1=f_{1,1}|f_{1,1}|d\mu_1=(f_{1,1}\sqrt{\rho_1})|f_{1,1}\sqrt{\rho_1}|d\nu$ 
and $u_1=f_{1,2}|f_{1,2}|d\mu_2=(f_{1,2}\sqrt{\rho_2})|f_{1,2}\sqrt{\rho_2}|d\nu$ hold. 
From Lemma 1, $\nu$-a.e., $f_{1,1}\sqrt{\rho_1}=f_{1,2}\sqrt{\rho_2}$ holds. 
Similarly, $\nu$-a.e., $f_{2,1}\sqrt{\rho_1}=f_{2,2}\sqrt{\rho_2}$ holds. 
$v_1=(f_{1,1}+f_{2,1})|f_{1,1}+f_{2,1}|d\mu_1
=(f_{1,1}\sqrt{\rho_1}+f_{2,1}\sqrt{\rho_1})|f_{1,1}\sqrt{\rho_1}+f_{2,1}\sqrt{\rho_1}|d\nu
=(f_{1,2}\sqrt{\rho_2}+f_{2,2}\sqrt{\rho_2})|f_{1,2}\sqrt{\rho_2}+f_{2,2}\sqrt{\rho_2}|d\nu
=(f_{1,2}+f_{2,2})|f_{1,2}+f_{2,2}|d\mu_2=v_2$ holds. 
\hfill 
$\blacksquare$

{\bf Definition 5 (scalar multiple):} 
Let $\Omega$ be a measurable space. Let $c\in \mathbb C$ and $u\in CM(\Omega)$.
Then, there uniquely exists $v\in CM(\Omega)$ such that 
the following holds.  
There exist a $\sigma$-finite measure $\mu$ on $\Omega$ 
and $f\in L^2(\mu)$ such that 
$u=f|f|d\mu$ and $v=(cf)|cf|d\mu$ hold.
Let $cu \in CM(\Omega)$ be defined as $cu:=v$. 

{\bf Proof:} 
From Lemma 1, existence is easy. 
Similar to Definition 4, 
uniqueness can be confirmed. 
\hfill 
$\blacksquare$

{\bf Definition 6 (inner product complex measure):} 
Let $\Omega$ be a measurable space. Let $u_1, u_2 \in CM(\Omega)$.
Then, there uniquely exists $v\in CM(\Omega)$ such that 
the following holds.  
There exist a $\sigma$-finite measure $\mu$ on $\Omega$ 
and $f_1, f_2\in L^2(\mu)$ such that 
$u_1=f_1|f_1|d\mu$, $u_2=f_2|f_2|d\mu$ 
and $v=f_1\overline{f_2}d\mu$ hold.
Let $\langle\langle u_1,u_2\rangle\rangle\in CM(\Omega)$ 
be defined as $\langle\langle u_1,u_2\rangle\rangle:=v$. 

{\bf Proof:} 
Similar to Definition 4, 
existence and uniqueness can be confirmed. 
\hfill 
$\blacksquare$

{\bf Definition 7 (inner product):} 
Let $\Omega$ be a measurable space. Let $u_1, u_2 \in CM(\Omega)$.
Then, let $\langle u_1,u_2\rangle_{CM(\Omega)} \in \mathbb C$ 
denote 
$\langle\langle u_1,u_2\rangle\rangle(\Omega)$.  
\hfill 
$\square$

{\bf Proposition 8:} 
Let $\Omega$ be a measurable space. 
Then, $CM(\Omega)$ is Hilbert space. 

{\bf Proof:} 
Similar to Definition 4, it can be confirmed 
that $CM(\Omega)$ is an inner product space. 
Let $\{u_n\}_n$ be Cauchy sequence in $CM(\Omega)$. 
We show that there exists $v\in CM(\Omega)$ such that 
$\lim_n \|u_n-v\|_{CM(\Omega)}=0$ holds. From Lemma 1 and Lemma 2, 
there exist a finite measure $\nu$ and a sequence $\{f_n\}_n$ 
in $L^2(\nu)$ such that for any $n$, $u_n=f_n|f_n|d\nu$ holds. 
So, because of $\|f_k-f_l\|_{L^2(\nu)}=\|u_k-u_l\|_{CM(\Omega)}$, 
$\{f_n\}_n$ is Cauchy sequence in $L^2(\nu)$. 
There exists $g\in L^2(\nu)$ such that 
$\lim_n \|f_n-g\|_{L^2(\nu)}=0$ holds. 
Then, because of $\|u_n-g|g|d\nu\|_{CM(\Omega)}
=\|f_n-g\|_{L^2(\nu)}$, 
$\lim_n \|u_n-g|g|d\nu\|_{CM(\Omega)}=0$ holds. 
\hfill 
$\blacksquare$

{\bf Remark:} 
(1) Let $M$ be a $C^\infty$-manifold. 
H\"{o}rmander ([6]) defined densities of order $\frac{1}{p}$ on $M$ ($\S2.4$) 
and further defined an inner product for densities of order $\frac{1}{2}$ on $M$ ($\S4.2$). 
Let $(\cdot,\cdot)_M$ be the inner product defined by H\"{o}rmander. 
Then, $(\cdot,\cdot)_M$ can be understood as a special case of the inner product of $CM(M)$. 
For simplicity, assume that $M$ is compact. 
Let $\omega$ be a volume element of $M$. 
Then, the positive square root $\sqrt{\omega}$ is a $C^\infty$-density of order $\frac{1}{2}$.  
For any $C^\infty$-density $u$ of order $\frac{1}{2}$, 
there uniquely exists a complex-valued $C^\infty$-function $f$ 
such that $u=f\sqrt{\omega}$ holds. 
For any complex-valued $C^\infty$-functions $f_1$ and $f_2$, 
$(f_1\sqrt{\omega},f_2\sqrt{\omega})_M=\int_Mf_1\overline{f_2}d\omega
=\langle f_1|f_1|d\omega, f_2|f_2|d\omega\rangle_{CM(M)}$ holds. 

(2) Let $\Omega$ be a measurable space. 
As Ay, Jost, L\^{e} and Schwachh\"{o}fer (Remark 3.5 of [1]) noted, 
Neveu ([11]) indicated that some real Hilbert space $S^{\frac{1}{2}}(\Omega)$ 
is naturally defined. 
On the other hand, the set $RM(\Omega)$ of all real measures on $\Omega$ is a real Hilbert space. 
It can be seen 
that $S^{\frac{1}{2}}(\Omega)$ and $RM(\Omega)$ are naturally isomorphic as real Hilbert spaces. 
Therefore, not only is a statistical manifold on $\Omega$ a submanifold of $RM(\Omega)$, 
but the inner product of $RM(\Omega)$ induces Fisher-Rao metric. 
\hfill 
$\square$

\newpage

\sect{Self-adjoint operator $\sqrt{-1}\frac{\partial}{\partial x_k}$}
As $\mathbb R$ is the measurable space whose set of all measurable sets is the topological $\sigma$-algebra, 
$\mathbb R^\infty$ denotes the product measurable space $\prod_{n\in\mathbb N}\mathbb R$.  

{\bf Lemma 9:} 
Let $\alpha, \beta \in \mathbb R^\infty$. Let $E$ be a subset of $\mathbb R^\infty$. 
Let $E_\alpha:=\{x\in\mathbb R^\infty|x-\alpha\in E\}$ and 
$E_\beta:=\{x\in\mathbb R^\infty|x-\beta\in E\}$. 
Then, if $E_\alpha$ is a measurable set of $\mathbb R^\infty$, 
then $E_\beta$ is a measurable set of $\mathbb R^\infty$. 

{\bf Proof:} 
Let $\cal M$ be the set of all subsets $D$ of $\mathbb R^\infty$ such that 
$\{x\in\mathbb R^\infty|x-(\beta-\alpha)\in D\}$
is a measurable set of $\mathbb R^\infty$. Then, $\cal M$ is a $\sigma$-algebra on $\mathbb R^\infty$. 
Furthermore, for any $k\in \mathbb N$ and any measurable set $B$ of $\mathbb R$, 
$(\prod_{n\in\{k\}}B)\times(\prod_{n\in\mathbb N\setminus\{k\}}\mathbb R) \in \cal M$ holds. 
So, if $D$ is a measurable set of $\mathbb R^\infty$, then $D\in \cal M$ holds. 
In particular, $E_\alpha\in \cal M$ holds. 
Therefore, since on the other hand, $E_\beta=\{x\in\mathbb R^\infty|x-(\beta-\alpha)\in E_\alpha\}$ holds, 
$E_\beta$ is a measurable set of $\mathbb R^\infty$. 
\hfill 
$\blacksquare$

{\bf Lemma 10:} 
Let $a\in \mathbb R^\infty$ and $T\in CM(\mathbb R^\infty)$. 
Let $T_a$ be the map from the set of all subsets $E$ 
of $\mathbb R^\infty$ such that 
$\{x\in\mathbb R^\infty|x+a\in E\}$
is a measurable set of $\mathbb R^\infty$ to $\mathbb C$ defined as 
$$T_a(E):=T(\{x\in\mathbb R^\infty|x+a\in E\}).$$
Then, $T_a\in CM(\mathbb R^\infty)$ holds. 

{\bf Proof:} 
From Lemma 9, it is easy. 
\hfill 
$\blacksquare$

{\bf Definition 11 (translation):}
Let $a\in \mathbb R^\infty$. Then, let a map $\tau_a$ 
from $CM(\mathbb R^\infty)$ to $CM(\mathbb R^\infty)$ be defined as 
$$(\tau_aT)(E):=T(\{x\in\mathbb R^\infty|x+a\in E\}).$$
\hfill 
$\square$

{\bf Lemma 12:} 
Let $a\in \mathbb R^\infty$. Let $\mu$ be a measure on $\mathbb R^\infty$. 
Let $\mu_a$ be the map from the set of all subsets $E$ 
of $\mathbb R^\infty$ such that 
$\{x\in\mathbb R^\infty|x+a\in E\}$
is a measurable set of $\mathbb R^\infty$ to $[0,+\infty]$ defined as 
$$\mu_a(E):=\mu(\{x\in\mathbb R^\infty|x+a\in E\}).$$
Then, $\mu_a$ is a measure on $\mathbb R^\infty$. 
Let $f$ be a complex measurable function on $\mathbb R^\infty$. 
Let $f_a$ be the map from $\mathbb R^\infty$ to $\mathbb C$ defined as 
$$f_a(x):=f(x-a).$$
Then, $f_a$ is a complex measurable function on $\mathbb R^\infty$. 
If $f\in L^1(\mu)$ holds, then $f_a\in L^1(\mu_a)$ and $\tau_a(fd\mu)=f_ad\mu_a$ hold.

{\bf Proof:} 
Although it is a natural result, we include the proof just in case. 

From Lemma 9, $\mu_a$ is a measure and $f_a$ is a measurable function. 

We show that if $f$ is non-negative, then $\tau_a(fd\mu)=f_ad\mu_a$ holds. 
There exists a monotonic non-decreasing sequence $\{g_n\}_n$ 
of non-negative simple measurable functions 
such that for any $x$, $\lim_n |g_n(x)-f(x)|=0$ holds. 
For $n$, let $g_{a,n}$ be the map defined as 
$g_{a,n}(x):=g_n(x-a)$. Then,  $\{g_{a,n}\}_n$ is a monotonic non-decreasing sequence  
of non-negative simple measurable functions 
such that for any $x$, $\lim_n |g_{a,n}(x)-f_a(x)|=0$ holds. 
So, for any measurable set $E$, $(\tau_a(fd\mu))(E)
=(fd\mu)(\{x\in\mathbb R^\infty|x+a\in E\})
=\int_{\{x\in\mathbb R^\infty|x+a\in E\}}fd\mu
=\lim_n (\int_{\{x\in\mathbb R^\infty|x+a\in E\}}g_nd\mu)
=\lim_n (\int_{E} g_{a,n}d\mu_a)
=\int_{E} f_ad\mu_a$ holds. 
If $f$ is non-negative, then $\tau_a(fd\mu)=f_ad\mu_a$ holds. 

There exist $h_1, h_2, h_3, h_4\in L^1(\mu)$ such that $h_1, h_2, h_3$ and $h_4$ 
are non-negative and 
$$f=(h_1-h_2)+\sqrt{-1}(h_3-h_4)$$ 
holds. Let $h_{a,1}(x):=h_1(x-a), h_{a,2}(x):=h_2(x-a), h_{a,3}(x):=h_3(x-a)$ and $h_{a,4}(x):=h_4(x-a)$. 
Then, $\tau_a(h_1d\mu)=h_{a,1}d\mu_a$, 
$\tau_a(h_2d\mu)=h_{a,2}d\mu_a$, 
$\tau_a(h_3d\mu)=h_{a,3}d\mu_a$ and 
$\tau_a(h_4d\mu)=h_{a,4}d\mu_a$ hold. 
Furthermore, $$f_a=(h_{a,1}-h_{a,2})+\sqrt{-1}(h_{a,3}-h_{a,4})$$ holds. 
So, for any measurable set $E$, $(\tau_a(fd\mu))(E)
=(fd\mu)(\{x\in\mathbb R^\infty|x+a\in E\})
=\int_{\{x\in\mathbb R^\infty|x+a\in E\}}fd\mu
=(\int_{\{x\in\mathbb R^\infty|x+a\in E\}}h_1d\mu
-\int_{\{x\in\mathbb R^\infty|x+a\in E\}}h_2d\mu)
+\sqrt{-1}(\int_{\{x\in\mathbb R^\infty|x+a\in E\}}h_3d\mu
-\int_{\{x\in\mathbb R^\infty|x+a\in E\}}h_4d\mu)
=((h_1d\mu)(\{x\in\mathbb R^\infty|x+a\in E\})-(h_2d\mu)(\{x\in\mathbb R^\infty|x+a\in E\}))
+\sqrt{-1}((h_3d\mu)(\{x\in\mathbb R^\infty|x+a\in E\})-(h_4d\mu)(\{x\in\mathbb R^\infty|x+a\in E\}))
=((\tau_a(h_1d\mu))(E)-(\tau_a(h_2d\mu))(E))
+\sqrt{-1}((\tau_a(h_3d\mu))(E)-(\tau_a(h_4d\mu))(E))
=((h_{a,1}d\mu_a)(E)-(h_{a,2}d\mu_a)(E))
+\sqrt{-1}$


\noindent
$((h_{a,3}d\mu_a)(E)-(h_{a,4}d\mu_a)(E))
=(\int_Eh_{a,1}d\mu_a-\int_Eh_{a,2}d\mu_a)
+\sqrt{-1}(\int_Eh_{a,3}d\mu_a-\int_Eh_{a,4}d\mu_a)
=\int_Ef_ad\mu_a=(f_ad\mu_a)(E)$ holds. 
\hfill 
$\blacksquare$

{\bf Proposition 13:} 
Let $a\in \mathbb R^\infty$. 
Then, $\tau_a$ is a unitary operator in $CM(\mathbb R^\infty)$. 

{\bf Proof:} 
From Lemma 12, it is easy. 
\hfill 
$\blacksquare$

{\bf Definition 14 ($L^2_k(\mathbb R^\infty)$, $(H^1_k(\mathbb R^\infty), \frac{\partial}{\partial x_k})$):} 
Let $k\in \mathbb N$. Then, there uniquely exists $e_k\in\mathbb R^\infty$ 
such that $e_{k,k}=1$ holds and for any $n\in\mathbb N\setminus\{k\}$, $e_{k,n}=0$ holds. 

(1) Let $L^2_k(\mathbb R^\infty)$ denote the set of all $u\in CM(\mathbb R^\infty)$ 
such that $$\lim_{h\downarrow+0}\|\tau_{he_k}u-u\|_{CM(\mathbb R^\infty)}=0$$ holds. 

(2) For $u_1, u_2 \in L^2_k(\mathbb R^\infty)$, let $\langle u_1,u_2 \rangle_{L^2_k(\mathbb R^\infty)}\in \mathbb C$ 
be defined as $$\langle u_1,u_2 \rangle_{L^2_k(\mathbb R^\infty)}:=\langle u_1,u_2 \rangle_{CM(\mathbb R^\infty)}.$$ 

(3) Let $H^1_k(\mathbb R^\infty)$ denote the set of all $u\in CM(\mathbb R^\infty)$ 
such that the following holds. There uniquely exists $v\in CM(\mathbb R^\infty)$ such that 
$$\lim_{h\downarrow+0}\|\frac{\tau_{he_k}u-u}{h}-v\|_{CM(\mathbb R^\infty)}=0$$ holds. 

(4) Let $u\in H^1_k(\mathbb R^\infty)$. Then, there uniquely exists $v\in CM(\mathbb R^\infty)$ such that 
$$\lim_{h\downarrow+0}\|\frac{\tau_{he_k}u-u}{h}-v\|_{CM(\mathbb R^\infty)}=0$$ holds. 
Let $\frac{\partial u}{\partial x_k}\in CM(\mathbb R^\infty)$ be defined as $\frac{\partial u}{\partial x_k}:=-v$. 
\hfill 
$\square$

{\bf Lemma 15:} 
Let $k\in \mathbb N$. Then, the followings hold. 

(1) $H^1_k(\mathbb R^\infty)\subset L^2_k(\mathbb R^\infty)$ holds. 

(2) $L^2_k(\mathbb R^\infty)$ is the set of all $u\in CM(\mathbb R^\infty)$ 
such that the following holds. 
For any $\varepsilon>0$, there exists $\delta>0$ such that 
for any $t_0, t_1, t_2\in \mathbb R$, if $|t_1-t_2|<\delta$ holds, 
then $\|\tau_{t_1e_k}(\tau_{t_0e_k}u)-\tau_{t_2e_k}(\tau_{t_0e_k}u)\|_{CM(\mathbb R^\infty)}<\varepsilon$ holds. 

(3) $L^2_k(\mathbb R^\infty)$ is a closed linear subspace of $CM(\mathbb R^\infty)$. 

(4) Let $u\in H^1_k(\mathbb R^\infty)$. Then, $-\frac{\partial u}{\partial x_k}\in L^2_k(\mathbb R^\infty)$ holds. 

{\bf Proof:} 
(1) It is easy. 

(2) Because of $\tau_{t_1e_k}(\tau_{t_0e_k}u)-\tau_{t_2e_k}(\tau_{t_0e_k}u)
=-\tau_{t_0e_k}\tau_{t_1e_k}(\tau_{(t_2-t_1)e_k}u-u)
=+\tau_{t_0e_k}\tau_{t_2e_k}(\tau_{(t_1-t_2)e_k}u-u)$, from Proposition 13, it follows. 

(3) From Proposition 13, it is easy. 

(4) From (1), (2) and (3), it follows. 
\hfill 
$\blacksquare$

{\bf Theorem 16:} 
Let $k\in \mathbb N$. Then, the followings hold. 

(1) The domain of $\sqrt{-1}\frac{\partial}{\partial x_k}$ is $H^1_k(\mathbb R^\infty)$. 

(2) $\sqrt{-1}\frac{\partial}{\partial x_k}$ is a self-adjoint operator in $L^2_k(\mathbb R^\infty)$. 

{\bf Proof:} 
(1) It is obvious. 

(2) From Proposition 8, Proposition 13 and Lemma 15, it follows (by Stone theorem). 
\hfill 
$\blacksquare$

\newpage

\sect{Non-negative self-adjoint operator $-\triangle_{\mathbb R^\infty}$}
{\bf Definition 17 ($H^1(\mathbb R^\infty)$):} 
(1) Let $H^1(\mathbb R^\infty)$ denote the set of all $u\in\cap_{k\in\mathbb N}H^1_k(\mathbb R^\infty)$ 
such that 
$$\sum_{k\in\mathbb N}\langle\sqrt{-1}\frac{\partial u}{\partial x_k},
\sqrt{-1}\frac{\partial u}{\partial x_k}\rangle_{L^2_k(\mathbb R^\infty)}<+\infty$$
holds. 

(2) For $u_1, u_2\in H^1(\mathbb R^\infty)$, let $\langle u_1, u_2 \rangle_{H^1(\mathbb R^\infty)}\in\mathbb C$ 
be defined as
$$\langle u_1, u_2 \rangle_{H^1(\mathbb R^\infty)}
:=\langle u_1, u_2 \rangle_{CM(\mathbb R^\infty)}
+\sum_{k\in\mathbb N}
\langle\sqrt{-1}\frac{\partial u_1}{\partial x_k},
\sqrt{-1}\frac{\partial u_2}{\partial x_k}\rangle_{L^2_k(\mathbb R^\infty)}.$$
\hfill 
$\square$

{\bf Proposition 18:} 
$H^1(\mathbb R^\infty)$ is Hilbert space. 

{\bf Proof:} 
It is easy to see that $H^1(\mathbb R^\infty)$ is an inner product space. 
Let $\{u_n\}_n$ be Cauchy sequence in $H^1(\mathbb R^\infty)$. 
Then, we show that there exists $v\in H^1(\mathbb R^\infty)$ such that 
$\lim_n\|u_n-v\|_{H^1(\mathbb R^\infty)}=0$ holds. 

There exist $v$ and $\{w_k\}_{k\in\mathbb N}$ such that 
$\lim_n\|u_n-v\|_{CM(\mathbb R^\infty)}=0$ holds and for any $k\in\mathbb N$, 
$\lim_n\|\sqrt{-1}\frac{\partial u_n}{\partial x_k}-w_k\|_{L^2_k(\mathbb R^\infty)}=0$ holds. 
For any $k\in\mathbb N$, because $\{u_n\}_n$ is Cauchy sequence in $L^2_k(\mathbb R^\infty)$, 
$\lim_n\|u_n-v\|_{L^2_k(\mathbb R^\infty)}=0$ holds. So, from Theorem 16, 
(because self-adjoint operators are closed,) 
for any $k\in\mathbb N$, $\sqrt{-1}\frac{\partial v}{\partial x_k}=w_k$ holds. 
Hence, for any $k\in\mathbb N$, $\lim_n\|\sqrt{-1}\frac{\partial u_n}{\partial x_k}
-\sqrt{-1}\frac{\partial v}{\partial x_k}\|_{L^2_k(\mathbb R^\infty)}=0$ holds. 
Because for any $n, K \in\mathbb N$, 
$$\sum_{k=1}^K
\langle\sqrt{-1}\frac{\partial u_n}{\partial x_k},
\sqrt{-1}\frac{\partial u_n}{\partial x_k}\rangle_{L^2_k(\mathbb R^\infty)}
\leq \sup_{m\in\mathbb N}\langle u_m,u_m\rangle_{H^1(\mathbb R^\infty)}<+\infty$$
holds, for any $K \in\mathbb N$, 
$$\sum_{k=1}^K
\langle\sqrt{-1}\frac{\partial v}{\partial x_k},
\sqrt{-1}\frac{\partial v}{\partial x_k}\rangle_{L^2_k(\mathbb R^\infty)}
\leq \sup_{m\in\mathbb N}\langle u_m,u_m\rangle_{H^1(\mathbb R^\infty)}<+\infty$$
holds. Therefore, $v\in H^1(\mathbb R^\infty)$ holds. 

Let $\varepsilon>0$. 
Then, there exists $N\in \mathbb N$ such that 
$$n\geq N, m\geq N \ \ \ \Longrightarrow \ \ \ \|u_n-u_m\|_{H^1(\mathbb R^\infty)}<\frac{\varepsilon}{2}$$
holds. For any $K\in\mathbb N$, 
$$n\geq N, m\geq N$$
$$\Longrightarrow$$
$$\langle u_n-u_m, u_n-u_m \rangle_{CM(\mathbb R^\infty)}
+\sum_{k=1}^K
\langle\sqrt{-1}\frac{\partial (u_n-u_m)}{\partial x_k},
\sqrt{-1}\frac{\partial (u_n-u_m)}{\partial x_k}\rangle_{L^2_k(\mathbb R^\infty)}$$
$$<(\frac{\varepsilon}{2})^2$$
holds. So, because of
$$n\geq N$$
$$\Longrightarrow$$
$$\langle u_n-v, u_n-v \rangle_{CM(\mathbb R^\infty)}
+\sum_{k=1}^K
\langle\sqrt{-1}\frac{\partial (u_n-v)}{\partial x_k},
\sqrt{-1}\frac{\partial (u_n-v)}{\partial x_k}\rangle_{L^2_k(\mathbb R^\infty)}$$
$$\leq (\frac{\varepsilon}{2})^2,$$
$$n\geq N \ \ \ \Longrightarrow \ \ \ 
\|u_n-v\|_{H^1(\mathbb R^\infty)}\leq\frac{\varepsilon}{2}<\varepsilon$$
holds. 
\hfill 
$\blacksquare$

{\bf Definition 19 ($L^2(\mathbb R^\infty)$):} 
Let $L^2(\mathbb R^\infty)$ denote the closure of $H^1(\mathbb R^\infty)$ in $CM(\mathbb R^\infty)$. 
For $u_1, u_2 \in L^2(\mathbb R^\infty)$, let $\langle u_1, u_2 \rangle_{L^2(\mathbb R^\infty)}\in\mathbb C$ 
be defined as $\langle u_1, u_2 \rangle_{L^2(\mathbb R^\infty)}
:=\langle u_1, u_2 \rangle_{CM(\mathbb R^\infty)}$. 
\hfill 
$\square$

{\bf Proposition 20:} 
$L^2(\mathbb R^\infty)$ is a closed linear subspace of $CM(\mathbb R^\infty)$. 
$L^2(\mathbb R^\infty)$ is Hilbert space. 
$H^1(\mathbb R^\infty)$ is a dense linear subspace of $L^2(\mathbb R^\infty)$. 
For any $u\in H^1(\mathbb R^\infty)$, 
$\|u\|_{L^2(\mathbb R^\infty)}\leq \|u\|_{H^1(\mathbb R^\infty)}$ holds. 

{\bf Proof:} 
It is easy. 
\hfill 
$\blacksquare$

{\bf Definition 21 (elliptic equation, $(1-\triangle_{\mathbb R^\infty})^{-1}$):} 
Let $f\in L^2(\mathbb R^\infty)$. Then, there uniquely exists $u\in H^1(\mathbb R^\infty)$ such that 
for any $\varphi\in H^1(\mathbb R^\infty)$, 
$$\langle u, \varphi \rangle_{H^1(\mathbb R^\infty)}
=\langle f, \varphi \rangle_{L^2(\mathbb R^\infty)}$$ 
holds. Let $(1-\triangle_{\mathbb R^\infty})^{-1}f\in H^1(\mathbb R^\infty)$ be defined as 
$(1-\triangle_{\mathbb R^\infty})^{-1}f:=u$. 

{\bf Proof:} 
From Proposition 18 and Proposition 20,  
existence and uniqueness are easy (by Riesz theorem). 
\hfill 
$\blacksquare$

{\bf Proposition 22:} 
For any $f\in L^2(\mathbb R^\infty)$, 
$\|(1-\triangle_{\mathbb R^\infty})^{-1}f\|_{H^1(\mathbb R^\infty)}
\leq \|f\|_{L^2(\mathbb R^\infty)}$ holds. 
$(1-\triangle_{\mathbb R^\infty})^{-1}$ 
is an injective self-adjoint operator in $L^2(\mathbb R^\infty)$.

{\bf Proof:} 
$\langle f_1, (1-\triangle_{\mathbb R^\infty})^{-1}f_2 \rangle_{L^2(\mathbb R^\infty)}
=\langle (1-\triangle_{\mathbb R^\infty})^{-1}f_1, 
(1-\triangle_{\mathbb R^\infty})^{-1}f_2 \rangle_{H^1(\mathbb R^\infty)}
=\overline{\langle (1-\triangle_{\mathbb R^\infty})^{-1}f_2, 
(1-\triangle_{\mathbb R^\infty})^{-1}f_1 \rangle_{H^1(\mathbb R^\infty)}}
=\overline{\langle f_2, (1-\triangle_{\mathbb R^\infty})^{-1}f_1 \rangle_{L^2(\mathbb R^\infty)}}
=\langle (1-\triangle_{\mathbb R^\infty})^{-1}f_1, f_2 \rangle_{L^2(\mathbb R^\infty)}$
holds. 

$\|(1-\triangle_{\mathbb R^\infty})^{-1}f\|_{H^1(\mathbb R^\infty)}^2
=\langle f, (1-\triangle_{\mathbb R^\infty})^{-1}f \rangle_{L^2(\mathbb R^\infty)}
\leq \|f\|_{L^2(\mathbb R^\infty)}\|(1-\triangle_{\mathbb R^\infty})^{-1}f\|_{L^2(\mathbb R^\infty)}
$ holds. 

There exists a sequence $\{g_n\}_n$ in $H^1(\mathbb R^\infty)$ such that 
$\lim_n\|g_n-f\|_{L^2(\mathbb R^\infty)}=0$ holds. 
So, if $(1-\triangle_{\mathbb R^\infty})^{-1}f=0$ holds, then 
$\langle f, f\rangle_{L^2(\mathbb R^\infty)}
=\lim_n \langle f, g_n\rangle_{L^2(\mathbb R^\infty)}
=\lim_n \langle (1-\triangle_{\mathbb R^\infty})^{-1}f, g_n\rangle_{H^1(\mathbb R^\infty)}=0$ holds. 
\hfill 
$\blacksquare$

{\bf Definition 23 ($\triangle_{\mathbb R^\infty}$):} 
Let $u$ be an element of the range of $(1-\triangle_{\mathbb R^\infty})^{-1}$. 
Then, let $\triangle_{\mathbb R^\infty}u\in L^2(\mathbb R^\infty)$ be defined as 
$$\triangle_{\mathbb R^\infty}u:=u-((1-\triangle_{\mathbb R^\infty})^{-1})^{-1}u.$$
\hfill 
$\square$

{\bf Theorem 24:} 
(1) The domain of $-\triangle_{\mathbb R^\infty}$ is a linear subspace of $H^1(\mathbb R^\infty)$. 
(2) $-\triangle_{\mathbb R^\infty}$ is a non-negative self-adjoint operator in $L^2(\mathbb R^\infty)$. 

{\bf Proof:} 
(1) The range of $(1-\triangle_{\mathbb R^\infty})^{-1}$ is a linear subspace of $H^1(\mathbb R^\infty)$. 

(2) From Proposition 22, $-\triangle_{\mathbb R^\infty}$ is self-adjoint. 
$\langle -\triangle_{\mathbb R^\infty}u, u\rangle_{L^2(\mathbb R^\infty)}
=\langle ((1-\triangle_{\mathbb R^\infty})^{-1})^{-1}u, u\rangle_{L^2(\mathbb R^\infty)}
-\langle u, u\rangle_{L^2(\mathbb R^\infty)}
=\langle u, u\rangle_{H^1(\mathbb R^\infty)}
-\langle u, u\rangle_{L^2(\mathbb R^\infty)}\geq 0$ holds. 
\hfill 
$\blacksquare$

\newpage

\sect{Product complex measure} 
We may use the following simple lemma without specific mention. 

{\bf Lemma 25:} 
Let $\Omega_1$ and $\Omega_2$ be measurable spaces. 
Let $\mu_1$ be a $\sigma$-finite measure on $\Omega_1$. 
Let $\mu_2$ be a $\sigma$-finite measure on $\Omega_2$. 
Let $\rho_1$ be a $[0,+\infty)$-valued measurable function on $\Omega_1$. 
Let $\rho_2$ be a $[0,+\infty)$-valued measurable function on $\Omega_2$. 
Then, $\rho_1\rho_2d\mu_1d\mu_2=(\rho_1d\mu_1)(\rho_2d\mu_2)$ holds. 

{\bf Proof:} 
Although it is a natural result, we include the proof just in case. 
Let $E_1$ be a measurable set of $\Omega_1$. 
Let $E_2$ be a measurable set of $\Omega_2$. 
Then, $(\rho_1\rho_2d\mu_1d\mu_2)(E_1\times E_2)
=\int_{E_1\times E_2}\rho_1\rho_2d\mu_1d\mu_2
=\int_{E_2}(\int_{E_1}\rho_1\rho_2d\mu_1)d\mu_2
=(\int_{E_1}\rho_1d\mu_1)(\int_{E_2}\rho_2d\mu_2)
=((\rho_1d\mu_1)(E_1))((\rho_2d\mu_2)(E_2))$ 
holds. 
\hfill 
$\blacksquare$

{\bf Definition 26 (product complex measure):} 
Let $u_1$ and $u_2$ be complex measures. 
Then, there uniquely exists a complex measure $v$ 
such that the following holds. 
There exist $\sigma$-finite measures $\mu_1, \mu_2$ 
and $f_1\in L^1(\mu_1), f_2\in L^1(\mu_2)$ such that 
$u_1=f_1d\mu_1$, $u_2=f_2d\mu_2$ 
and $v=f_1f_2d\mu_1d\mu_2$ hold. 
Let the complex measure $u_1\cdot u_2$ be defined as $u_1\cdot u_2:=v$. 

{\bf Proof:} 
Existence is easy. We show uniqueness. 
Suppose that $(v_1,\mu_{1,1},\mu_{2,1},$


\noindent 
$f_{1,1},f_{2,1})$ and 
$(v_2,\mu_{1,2},\mu_{2,2},f_{1,2},f_{2,2})$ each satisfy the condition of the definition. 
Then, there exist finite measures $\nu_1, \nu_2$ 
and $[0,+\infty)$-valued measurable functions $\rho_{1,1}, \rho_{2,1}, \rho_{1,2}, \rho_{2,2}$ such that 
$\mu_{1,1}=\rho_{1,1}d\nu_1$, $\mu_{1,2}=\rho_{1,2}d\nu_1$, 
$\mu_{2,1}=\rho_{2,1}d\nu_2$ and $\mu_{2,2}=\rho_{2,2}d\nu_2$ hold.  
So, because $f_{1,1}\rho_{1,1}d\nu_1=u_1=f_{1,2}\rho_{1,2}d\nu_1$ and  
$f_{2,1}\rho_{2,1}d\nu_2=u_2=f_{2,2}\rho_{2,2}d\nu_2$ hold,  
$d\nu_1d\nu_2$-a.e., $f_{1,1}\rho_{1,1}f_{2,1}\rho_{2,1}
=f_{1,2}\rho_{1,2}f_{2,2}\rho_{2,2}$ holds. 
From Lemma 25, $v_1=f_{1,1}f_{2,1}d\mu_{1,1}d\mu_{2,1}
=f_{1,1}\rho_{1,1}f_{2,1}\rho_{2,1}d\nu_1d\nu_2
=f_{1,2}\rho_{1,2}f_{2,2}\rho_{2,2}d\nu_1d\nu_2
=f_{1,2}f_{2,2}d\mu_{1,2}d\mu_{2,2}
=v_2$ holds. 
\hfill 
$\blacksquare$

{\bf Proposition 27:} 
Let $u_1, u_2$ and $u_3$ be complex measures. 
Then, $$(u_1\cdot u_2)\cdot u_3=u_1\cdot(u_2\cdot u_3)$$ holds. 

{\bf Proof:} 
It is easy. 
\hfill 
$\blacksquare$

{\bf Proposition 28:} 
Let $\Omega_1$ and $\Omega_2$ be measurable spaces. 
Then, the followings hold. 

(1) Let $u_1, v_1\in CM(\Omega_1)$ and $u_2, v_2\in CM(\Omega_2)$. Then, 
$$\langle\langle u_1\cdot u_2, v_1\cdot v_2\rangle\rangle=\langle\langle u_1, v_1\rangle\rangle
\cdot\langle\langle u_2, v_2\rangle\rangle,$$
$$\langle u_1\cdot u_2, v_1\cdot v_2\rangle_{CM(\Omega_1\times\Omega_2)}
=\langle u_1, v_1\rangle_{CM(\Omega_1)}
\langle u_2, v_2\rangle_{CM(\Omega_2)}$$
hold. 

(2) The map $(u_1,u_2) \mapsto u_1\cdot u_2$ from $CM(\Omega_1)\times CM(\Omega_2)$ 
to $CM(\Omega_1\times\Omega_2)$ is bi-linear.

{\bf Proof:} 
(1) There exist finite measures $\mu_1, \mu_2$ and $f_1, g_1\in L^2(\mu_1), f_2, g_2\in L^2(\mu_2)$ such that 
$u_1=f_1|f_1|d\mu_1, v_1=g_1|g_1|d\mu_1, u_2=f_2|f_2|d\mu_2$ and $v_2=g_2|g_2|d\mu_2$ hold. 
Then, from 
$\langle\langle u_1, v_1\rangle\rangle
=f_1\overline{g_1}d\mu_1$, 
$\langle\langle u_2, v_2\rangle\rangle
=f_2\overline{g_2}d\mu_2$ and 
$\langle\langle u_1\cdot u_2, v_1\cdot v_2\rangle\rangle
=f_1\overline{g_1}f_2\overline{g_2}d\mu_1d\mu_2$, it follows.  

(2) It is easy. 
\hfill 
$\blacksquare$

We may use (2) of the following remark without specific mention, 
as it is a natural result. 

{\bf Remark:} 
(1) Let $\Omega$ be a set. Let $\cal A$ be an algebra on $\Omega$. 
Let $\cal B$ be the smallest $\sigma$-algebra on $\Omega$ 
such that $\cal A\subset \cal B$ holds. Let $u_1$ and $u_2$ 
be real measures on the measurable space 
whose set of all measurable sets is $\cal B$. 
Suppose that for any $E\in \cal A$, $u_1(E)=u_2(E)$ holds.  
Then, $u_1=u_2$ holds.  

(2) Let $\Omega_1$ and $\Omega_2$ be measurable spaces. 
Let $u_1, u_2\in CM(\Omega_1\times\Omega_2)$. 
Suppose that for any measurable set $E_1$ of $\Omega_1$ 
and any measurable set $E_2$ of $\Omega_2$, 
$u_1(E_1\times E_2)=u_2(E_1\times E_2)$ holds. 
Then, $u_1=u_2$ holds.

{\bf Proof:} 
(1) Let $|u_1|$ be the total variation of $u_1$. 
Let $|u_2|$ be the total variation of $u_2$. 
For a measurable set $E$, let $\mu(E):=|u_1|(E)+|u_2|(E)$. 
Then, there exist $[-1,+1]$-valued measurable functions $f_1$ and $f_2$ 
such that $u_1=f_1d\mu$ and $u_2=f_2d\mu$ hold. 
Let $E_{f_1<f_2}:=\{x\in\Omega | f_1(x)<f_2(x)\}$
and $E_{f_2<f_1}:=\{x\in\Omega | f_2(x)<f_1(x)\}$. 

Let $\varepsilon>0$. Then,  
from basic facts about Carath\'{e}odory outer measure 
(especially, Hopf extension theorem), there exists a sequence $\{E_n\}_{n=1}^\infty$ in $\cal A$ 
such that $E_{f_1<f_2}\subset \cup_n E_n$ and $\sum_n\mu(E_n)<\mu(E_{f_1<f_2})+\frac{\varepsilon}{2}$ hold. 
Let $E^\prime_n:=E_n\setminus \cup_{k=1}^{n-1}E_k$. 
Then, because 
$\sum_n\mu(E^\prime_n\setminus E_{f_1<f_2})
=\mu(\cup_n(E^\prime_n\setminus E_{f_1<f_2}))
=\mu((\cup_n E_n)\setminus E_{f_1<f_2})
=\mu(\cup_n E_n)-\mu(E_{f_1<f_2})<\frac{\varepsilon}{2}$ 
and $E^\prime_n\in\cal A$ hold, 
$\int_{E_{f_1<f_2}}|f_1-f_2|d\mu=\sum_n \int_{E^\prime_n\cap E_{f_1<f_2}}(f_2-f_1)d\mu
=\sum_n (\int_{E^\prime_n}(f_2-f_1)d\mu
-\int_{E^\prime_n\setminus E_{f_1<f_2}}(f_2-f_1)d\mu)
=\sum_n \int_{E^\prime_n\setminus E_{f_1<f_2}}(f_1-f_2)d\mu
\leq \sum_n \int_{E^\prime_n\setminus E_{f_1<f_2}}2d\mu<\varepsilon$ holds. 
Therefore, $\int_{E_{f_1<f_2}}|f_1-f_2|d\mu=0$ holds. Similarly, 
$\int_{E_{f_2<f_1}}|f_1-f_2|d\mu=0$ holds. 

(2) Let $\cal P$ be the set of all product sets 
of a measurable set of $\Omega_1$ 
and a measurable set of $\Omega_2$. 
Let $\cal A$ be the smallest algebra on $\Omega_1\times\Omega_2$ 
such that $\cal P \subset \cal A$ holds. 
Then, $\cal A$ is the set of all finite disjoint unions 
of elements of $\cal P$. So, for any $E\in \cal A$, $u_1(E)=u_2(E)$ holds. 
Let $\cal B$ be the smallest $\sigma$-algebra on $\Omega_1\times\Omega_2$ 
such that $\cal A\subset \cal B$ holds. 
Then, $\cal B$ is the set of all measurable sets of $\Omega_1\times\Omega_2$.
Therefore, from (1),  $\mathrm{Re}(u_1)=\mathrm{Re}(u_2)$ 
and $\mathrm{Im}(u_1)=\mathrm{Im}(u_2)$ hold. 

\hfill 
$\blacksquare$

\newpage

\sect{Embedding finite-dimensional evolution \, \, 
$\{e^{\triangle_{\mathbb R^N}t}\}_{t\in[0,+\infty)}$ and 
$\{e^{\sqrt{-1}\triangle_{\mathbb R^N}t}\}_{t\in(-\infty,+\infty)}$}
Let $N\in\mathbb N$. 
As $\mathbb R$ is the measurable space whose set of all measurable sets is the topological $\sigma$-algebra, 
let $\mathbb R^N$ denote the product measurable space $\prod_{n=1}^N\mathbb R$.  
Let $dx$ denote the (ordinary) measure $\prod_{n=1}^Ndx_n$ on $\mathbb R^N$. 
Let $\triangle_{\mathbb R^N}$ denote the (ordinary) Laplacian in $L^2(\mathbb R^N)$. 

In Section 3, we defined $\{\frac{\partial}{\partial x_k}\}_{k\in\mathbb N}$. 
However, in $\mathbb R^N\times\mathbb R^\infty$, the variables $x_1, x_2, \cdots, x_N$ conflict. 
For notational consistency, we introduce the following definition.

{\bf Definition 29 ($N$-shift):} 
Let $u\in CM(\mathbb R^\infty)$. 
Then, let $u^+_N$ denote the map from the set of all subsets $E$ 
of $\prod_{n=N+1}^\infty\mathbb R$ such that 
$\{\{x_n\}_{n\in\mathbb N}\in\mathbb R^\infty|\{x_{n-N}\}_{n=N+1}^\infty\in E\}$
is a measurable set of $\mathbb R^\infty$ to $\mathbb C$ defined as 
$$u^+_N(E):=u(\{\{x_n\}_{n\in\mathbb N}\in\mathbb R^\infty|\{x_{n-N}\}_{n=N+1}^\infty\in E\}).$$
\hfill 
$\square$ 

{\bf Lemma 30:} 
For any $u\in CM(\mathbb R^\infty)$, $u^+_N\in CM(\prod_{n=N+1}^\infty\mathbb R)$ holds. 
The map $u\mapsto u^+_N$ from $CM(\mathbb R^\infty)$ to $CM(\prod_{n=N+1}^\infty\mathbb R)$ 
is a unitary operator. 

{\bf Proof:} 
It is easy. 
\hfill 
$\blacksquare$

{\bf Definition 31 ($\otimes$):} 
Let $f\in L^2(\mathbb R^N)$ and $u\in CM(\mathbb R^\infty)$. 
Then, let $f\otimes u\in CM(\mathbb R^\infty)$ defined as 
$$f\otimes u:=(f|f|dx)\cdot(u^+_N).$$
\hfill 
$\square$ 

{\bf Lemma 32:} 
Let $f_1, f_2\in L^2(\mathbb R^N)$ and $u_1, u_2\in CM(\mathbb R^\infty)$. 
Then, 
$$\langle f_1\otimes u_1, f_2\otimes u_2\rangle_{CM(\mathbb R^\infty)}
=\langle f_1, f_2\rangle_{L^2(\mathbb R^N)}
\langle u_1, u_2\rangle_{CM(\mathbb R^\infty)}$$
holds. 

{\bf Proof:} 
From Proposition 28 (1) and Lemma 30, it is easy. 
\hfill 
$\blacksquare$

{\bf Lemma 33:} 
The map $(f,u) \mapsto f\otimes u$ from $L^2(\mathbb R^N)\times CM(\mathbb R^\infty)$ 
to $CM(\mathbb R^\infty)$ is bi-linear. 

{\bf Proof:} 
From Proposition 28 (2) and Lemma 30, it is easy. 
\hfill 
$\blacksquare$ 

{\bf Lemma 34:} 
Let $a\in \mathbb R^N$. 
Let $f$ be a complex measurable function on $\mathbb R^N$. 
Let $f_a$ be the map from $\mathbb R^N$ to $\mathbb C$ defined as 
$$f_a(x):=f(x-a).$$ 
Then, $f_a$ is a complex measurable function on $\mathbb R^N$. 
Let $E$ be a measurable set of $\mathbb R^N$. 
Let $$E_a:=\{x\in\mathbb R^N|x-a\in E\}.$$
Then, $E_a$ is a measurable set of $\mathbb R^N$.  
If $f\in L^1(\mathbb R^N)$ holds, then $f_a\in L^1(\mathbb R^N)$ 
and $\int_{E_a}f_adx=\int_{E}fdx$ hold.

{\bf Proof:} 
It is well known. 
\hfill 
$\blacksquare$ 

{\bf Lemma 35:} 
Let $\{a_n\}_{n\in\mathbb N}\in \mathbb R^\infty$. 
Let $f\in L^2(\mathbb R^N)$ and $u\in CM(\mathbb R^\infty)$. 
Let  $f_{\{a_n\}_{n=1}^N}$ be the map from $\mathbb R^N$ to $\mathbb C$ defined as 
$$f_{\{a_n\}_{n=1}^N}(x):=f(x-{\{a_n\}_{n=1}^N}).$$ 
Then, 
$$\tau_{\{a_n\}_{n\in\mathbb N}}(f\otimes u)
=f_{\{a_n\}_{n=1}^N}\otimes (\tau_{\{a_{N+n}\}_{n\in\mathbb N}}u)$$
holds.

{\bf Proof:} 
Although it is a natural result, we include the proof just in case. 

Let $F$ be a measurable set of $\mathbb R^N$. Let $G$ be a measurable set of $\prod_{n=N+1}^\infty\mathbb R$. 
Then, from $\{x\in\mathbb R^\infty|x+\{a_n\}_{n=1}^\infty\in F\times G\}
=\{x\in\mathbb R^N|x+\{a_n\}_{n=1}^N\in F\}\times 
\{x\in\mathbb \prod_{n=N+1}^\infty\mathbb R|x+\{a_n\}_{n=N+1}^\infty\in G\}$
and Lemma 34, 
$(\tau_{\{a_n\}_{n=1}^\infty}(f\otimes u))(F\times G)
=(f\otimes u)(\{x\in\mathbb R^\infty|x+\{a_n\}_{n=1}^\infty\in F\times G\})
=(\int_{x+\{a_n\}_{n=1}^N\in F}f|f|dx)
(u^+_N(\{x\in\mathbb \prod_{n=N+1}^\infty\mathbb R|x+\{a_n\}_{n=N+1}^\infty\in G\}))
=(\int_F f_{\{a_n\}_{n=1}^N}|f_{\{a_n\}_{n=1}^N}|dx)
((\tau_{\{a_{N+n}\}_{n=1}^\infty}u)^+_N(G))$
holds. So, from Remark (2) at the end of Section 5, 
$\tau_{\{a_n\}_{n=1}^\infty}(f\otimes u)=(f_{\{a_n\}_{n=1}^N}|f_{\{a_n\}_{n=1}^N}|dx)
\cdot((\tau_{\{a_{N+n}\}_{n=1}^\infty}u)^+_N)$ holds. 
\hfill 
$\blacksquare$

{\bf Lemma 36:} 
Let $k\in \mathbb N$, $u\in H^1_{k}(\mathbb R^\infty)$ and $f\in L^2(\mathbb R^N)$. 
Then, $f\otimes u\in H^1_{N+k}(\mathbb R^\infty)$ and 
$$\frac{\partial (f\otimes u)}{\partial x_{N+k}}=f\otimes \frac{\partial u}{\partial x_{k}}$$
hold.

{\bf Proof:} 
There uniquely exists $e_{k}\in\mathbb R^\infty$ 
such that $e_{k,k}=1$ holds and for any $n\in\mathbb N\setminus\{k\}$, $e_{k,n}=0$ holds. 
There uniquely exists $e_{N+k}\in\mathbb R^\infty$ 
such that $e_{N+k,N+k}=1$ holds and for any $n\in\mathbb N\setminus\{N+k\}$, $e_{N+k,n}=0$ holds. 
Then, from Lemma 35, Lemma 33 and Lemma 32, 
$\|\frac{\tau_{he_{N+k}}(f\otimes u)-f\otimes u}{h}+f\otimes \frac{\partial u}{\partial x_{k}}\|_{CM(\mathbb R^\infty)}
=\|\frac{f\otimes(\tau_{he_{k}}u)-f\otimes u}{h}+f\otimes \frac{\partial u}{\partial x_{k}}\|_{CM(\mathbb R^\infty)}
=\|f\otimes(\frac{\tau_{he_{k}}u-u}{h}+\frac{\partial u}{\partial x_{k}})\|_{CM(\mathbb R^\infty)}
=\|f\|_{L^2(\mathbb R^N)}\|\frac{\tau_{he_{k}}u-u}{h}+\frac{\partial u}{\partial x_{k}}\|_{CM(\mathbb R^\infty)}$ 
holds. 
\hfill 
$\blacksquare$

{\bf Lemma 37:} 
Let $k\in \{1,2,\cdots,N\}$. Then, there uniquely exists $e^N_k\in \mathbb R^N$ such that  
$e^N_{k,k}=1$ holds and for any $n\in \{1,2,\cdots,N\}\setminus\{k\}$, $e^N_{k,n}=0$ holds. 
Let $f,g \in L^2(\mathbb R^N)$.  Suppose that 
$$\lim_{h\downarrow+0}\|\frac{f(x-he^N_k)-f(x)}{h}+g(x)\|_{L^2(\mathbb R^N)}=0$$
holds. Let $u\in CM(\mathbb R^\infty)$. Then, $f\otimes u\in H^1_{k}(\mathbb R^\infty)$ and 
$$\frac{\partial (f\otimes u)}{\partial x_{k}}=g\otimes u$$
hold. 

{\bf Proof:} 
There uniquely exists $e_k\in \mathbb R^\infty$ such that  
$e_{k,k}=1$ holds and for any $n\in\mathbb N\setminus\{k\}$, $e_{k,n}=0$ holds. 
For $a\in \mathbb R^N$, let $f_{a}(x):=f(x-a)$. 
Then, from Lemma 35, Lemma 33 and Lemma 32, 
$\|\frac{\tau_{he_{k}}(f\otimes u)-f\otimes u}{h}+g\otimes u\|_{CM(\mathbb R^\infty)}
=\|\frac{f_{he^N_k}\otimes u-f\otimes u}{h}+g\otimes u\|_{CM(\mathbb R^\infty)}
=\|(\frac{f_{he^N_k}-f}{h}+g)\otimes u\|_{CM(\mathbb R^\infty)}
=\|\frac{f_{he^N_k}-f}{h}+g\|_{L^2(\mathbb R^N)}\|u\|_{CM(\mathbb R^\infty)}$ 
holds. 
\hfill 
$\blacksquare$

{\bf Lemma 38:}
(1) Let $f\in H^1(\mathbb R^N)$ and $u\in H^1(\mathbb R^\infty)$. 
Then, $f\otimes u\in H^1(\mathbb R^\infty)$ holds. 

(2) Let $f\in L^2(\mathbb R^N)$ and $u\in L^2(\mathbb R^\infty)$. 
Then, $f\otimes u\in L^2(\mathbb R^\infty)$ holds. 

{\bf Proof:} 
(1) From Lemma 37 and Lemma 36, $f\otimes u\in \cap_{k\in\mathbb N}H^1_k(\mathbb R^\infty)$ holds 
and for any $k\in\{N+1,N+2,\cdots\}$, 
$\frac{\partial (f\otimes u)}{\partial x_{k}}=f\otimes \frac{\partial u}{\partial x_{k-N}}$ 
holds. So, from Lemma 32, 
$\sum_{k=N+1}^\infty\|\frac{\partial (f\otimes u)}{\partial x_{k}}\|_{CM(\mathbb R^\infty)}^2
=\|f\|_{L^2(\mathbb R^N)}^2(\sum_{k\in\mathbb N}\|\frac{\partial u}{\partial x_{k}}\|_{CM(\mathbb R^\infty)}^2)
\leq \|f\|_{L^2(\mathbb R^N)}^2$


\noindent 
$\|u\|_{H^1(\mathbb R^\infty)}^2$ 
holds. 

(2) There exists a sequence $\{g_n\}_n$ in $H^1(\mathbb R^N)$ 
such that $\lim_n\|g_n-f\|_{L^2(\mathbb R^N)}=0$ holds. 
There exists a sequence $\{v_n\}_n$ in $H^1(\mathbb R^\infty)$ 
such that $\lim_n\|v_n-u\|_{L^2(\mathbb R^\infty)}=0$ holds. 
Then, from (1), for any $n$, $g_n\otimes v_n\in H^1(\mathbb R^\infty)$ holds. 
On the other hand, from Lemma 33 and Lemma 32, 
$\lim_n\|g_n\otimes v_n-f\otimes u\|_{CM(\mathbb R^\infty)}=0$ holds. 
\hfill 
$\blacksquare$

In order to examine $(\triangle_{\mathbb R^N}f)\otimes u$ and 
$f\otimes (\triangle_{\mathbb R^\infty}u)$, we introduce $\otimes$-contraction 
($\diamond_{N}$ and $\diamond_{\infty}$).

{\bf Definition 39 ($\diamond_{N}$, $\diamond_{\infty}$):} 
Let $v\in CM(\mathbb R^\infty)$. 

(1) Let $f\in L^2(\mathbb R^N)$. Then, there uniquely exists $w\in CM(\mathbb R^\infty)$ 
such that for any $u\in CM(\mathbb R^\infty)$, 
$$\langle f\otimes u, v\rangle_{CM(\mathbb R^\infty)}
=\langle u, w\rangle_{CM(\mathbb R^\infty)}$$
holds. Let $f \diamond_{N}v\in CM(\mathbb R^\infty)$ be defined as $f\diamond_{N}v:=w$.

(2) Let $u\in CM(\mathbb R^\infty)$. Then, there uniquely exists $g\in L^2(\mathbb R^N)$ 
such that for any $f\in L^2(\mathbb R^N)$, 
$$\langle f\otimes u, v\rangle_{CM(\mathbb R^\infty)}
=\langle f, g\rangle_{L^2(\mathbb R^N)}$$
holds. Let $u\diamond_{\infty}v\in L^2(\mathbb R^N)$ be defined as $u\diamond_{\infty}v:=g$.

{\bf Proof:} 
Existence and uniqueness are easy (by Riesz Theorem).
\hfill 
$\blacksquare$

{\bf Lemma 40:} 
(1) Let $f\in L^2(\mathbb R^N)$ and $v\in CM(\mathbb R^\infty)$. 
Then, $\|f\diamond_{N}v\|_{CM(\mathbb R^\infty)}
\leq \|f\|_{L^2(\mathbb R^N)}\|v\|_{CM(\mathbb R^\infty)}$ 
holds. 

(2) Let $u, v\in CM(\mathbb R^\infty)$. 
Then, $\|u\diamond_{\infty}v\|_{L^2(\mathbb R^N)}
\leq \|u\|_{CM(\mathbb R^\infty)}\|v\|_{CM(\mathbb R^\infty)}$ 
holds. 

{\bf Proof:} 
It is easy.
\hfill 
$\blacksquare$

{\bf Proposition 41:} 
(1) Let $k\in\mathbb N$, $v\in H^1_{N+k}(\mathbb R^\infty)$ 
and $f\in L^2(\mathbb R^N)$. Then, 
$f\diamond_{N}v\in H^1_{k}(\mathbb R^\infty)$ and 
$$\frac{\partial (f\diamond_{N}v)}{\partial x_k}
=f\diamond_{N}\frac{\partial v}{\partial x_{N+k}}$$
hold. 

(2) Let $k\in\mathbb \{1,2,\cdots,N\}$, $v\in H^1_{k}(\mathbb R^\infty)$ 
and $u\in CM(\mathbb R^\infty)$. 
Let $\frac{\partial}{\partial x_{N,k}}$ 
be the (ordinary) generalized partial differential operator in $L^2(\mathbb R^N)$. 
Then, $u\diamond_{\infty}v$ is an element of the domain of $\frac{\partial}{\partial x_{N,k}}$ and 
$$\frac{\partial (u\diamond_{\infty}v)}{\partial x_{N,k}}
=u\diamond_{\infty}\frac{\partial v}{\partial x_{k}}$$
holds.

{\bf Proof:} 
(1) From Lemma 36 and Theorem 16, because for any $u\in H^1_k(\mathbb R^\infty)$, 
$\langle u, f\diamond_{N}\frac{\partial v}{\partial x_{N+k}}\rangle_{CM(\mathbb R^\infty)}
=\langle f\otimes u, \frac{\partial v}{\partial x_{N+k}}\rangle_{CM(\mathbb R^\infty)}
=-\langle f\otimes \frac{\partial u}{\partial x_{k}}, v\rangle_{CM(\mathbb R^\infty)}
=-\langle \frac{\partial u}{\partial x_{k}}, f\diamond_{N}v\rangle_{CM(\mathbb R^\infty)}$ 
holds, $f\diamond_{N}\frac{\partial v}{\partial x_{N+k}}
=-(\frac{\partial}{\partial x_{k}})^*(f\diamond_{N}v)
=\frac{\partial (f\diamond_{N}v)}{\partial x_k}$ holds. 

(2) In the same way as (1), from Lemma 37 and Theorem 16, it follows.

\hfill 
$\blacksquare$

{\bf Lemma 42:} 
(1) Let $\varphi\in H^1(\mathbb R^\infty)$ and $f\in L^2(\mathbb R^N)$. Then, 
$f\diamond_{N}\varphi\in H^1(\mathbb R^\infty)$ holds. 

(2) Let $\varphi\in H^1(\mathbb R^\infty)$ and $u\in CM(\mathbb R^\infty)$. Then, 
$u\diamond_{\infty}\varphi\in H^1(\mathbb R^N)$ holds. 

(3) Let $f\in L^2(\mathbb R^N)$. Let $u$ be an element of the domain of $\triangle_{\mathbb R^\infty}$. 
Then, $f\otimes (\triangle_{\mathbb R^\infty}u)\in L^2(\mathbb R^\infty)$, 
$f\otimes u\in\cap_{k\in\mathbb N}H^1_{N+k}(\mathbb R^\infty)$ and  
$\sum_{k\in\mathbb N}\|\frac{\partial (f\otimes u)}{\partial x_{N+k}}\|_{L^2_{N+k}(\mathbb R^\infty)}^2<+\infty$ hold. 
For any $\varphi\in H^1(\mathbb R^\infty)$, 
$$-\langle f\otimes (\triangle_{\mathbb R^\infty}u), \varphi\rangle_{L^2(\mathbb R^\infty)}
=\sum_{k\in\mathbb N}\langle \frac{\partial (f\otimes u)}{\partial x_{N+k}}, 
\frac{\partial \varphi}{\partial x_{N+k}}\rangle_{L^2_{N+k}(\mathbb R^\infty)}$$
holds.

(4) Let $u\in L^2(\mathbb R^\infty)$ and $f\in H^2(\mathbb R^N)$. Then, 
$(\triangle_{\mathbb R^N}f)\otimes u\in L^2(\mathbb R^\infty)$ and  
$f\otimes u\in\cap_{k=1}^NH^1_{k}(\mathbb R^\infty)$ hold. 
For any $\varphi\in H^1(\mathbb R^\infty)$, 
$$-\langle (\triangle_{\mathbb R^N}f)\otimes u, \varphi\rangle_{L^2(\mathbb R^\infty)}
=\sum_{k=1}^N\langle \frac{\partial (f\otimes u)}{\partial x_{k}}, 
\frac{\partial \varphi}{\partial x_{k}}\rangle_{L^2_{k}(\mathbb R^\infty)}$$
holds. 

{\bf Proof:} 
(1) From Proposition 41 (1) and Lemma 40 (1), it follows. 

(2) From Proposition 41 (2), it follows. 

(3) From Lemma 38 (2), $f\otimes (\triangle_{\mathbb R^\infty}u)\in L^2(\mathbb R^\infty)$ holds. 
From Theorem 24 (1), $u\in H^1(\mathbb R^\infty)$ holds.  
So, from Lemma 36, $f\otimes u\in\cap_{k\in\mathbb N}H^1_{N+k}(\mathbb R^\infty)$ and  
$\sum_{k\in\mathbb N}\|\frac{\partial (f\otimes u)}{\partial x_{N+k}}\|_{L^2_{N+k}(\mathbb R^\infty)}^2<+\infty$ hold. 
From (1), Proposition 41 (1) and Lemma 36, 
$\langle f\otimes u, \varphi\rangle_{CM(\mathbb R^\infty)}
-\langle f\otimes (\triangle_{\mathbb R^\infty}u), \varphi\rangle_{L^2(\mathbb R^\infty)}
=\langle (1-\triangle_{\mathbb R^\infty})u, f\diamond_{N}\varphi\rangle_{L^2(\mathbb R^\infty)}
=\langle u, f\diamond_{N}\varphi\rangle_{H^1(\mathbb R^\infty)}
=\langle u, f\diamond_{N}\varphi\rangle_{CM(\mathbb R^\infty)}
+\sum_{k\in\mathbb N}\langle \frac{\partial u}{\partial x_k}, 
f\diamond_{N} \frac{\partial \varphi}{\partial x_{N+k}}\rangle_{CM(\mathbb R^\infty)}
=\langle f\otimes u, \varphi\rangle_{CM(\mathbb R^\infty)}
+\sum_{k\in\mathbb N}\langle \frac{\partial (f\otimes u)}{\partial x_{N+k}}, 
\frac{\partial \varphi}{\partial x_{N+k}}\rangle_{L^2_{N+k}(\mathbb R^\infty)}$
holds. 

(4) Similar to (3), from (2), Lemma 37, Lemma 38 (2) and Proposition 41 (2), it is shown. 
\hfill 
$\blacksquare$

{\bf Lemma 43:} 
Suppose that $-\infty<T_0<T_1<+\infty$ holds. Let $f\in C^1([T_0,T_1];L^2(\mathbb R^N))$
and $u\in C^1([T_0,T_1];CM(\mathbb R^\infty))$. 
Then, $f\otimes u\in C^1([T_0,T_1]$


\noindent 
$;CM(\mathbb R^\infty))$ and 
$$\frac{d}{dt}(f\otimes u)=(\frac{d}{dt}f)\otimes u+f\otimes(\frac{d}{dt}u)$$
hold. 

{\bf Proof:} 
From Lemma 32 and Lemma 33, it follows. 
\hfill 
$\blacksquare$

{\bf Theorem 44:} 
Let $f_0\in L^2(\mathbb R^N)$ and $u_0\in L^2(\mathbb R^\infty)$. 
Then, the followings hold. 

(1) Let $t\in [0,+\infty)$. Then,
$$e^{\triangle_{\mathbb R^\infty}t}(f_0\otimes u_0)
=(e^{\triangle_{\mathbb R^N}t}f_0)\otimes 
(e^{\triangle_{\mathbb R^\infty}t}u_0)$$
holds.

(2) Let $t\in (-\infty,+\infty)$. Then,
$$e^{\sqrt{-1}\triangle_{\mathbb R^\infty}t}(f_0\otimes u_0)
=(e^{\sqrt{-1}\triangle_{\mathbb R^N}t}f_0)\otimes 
(e^{\sqrt{-1}\triangle_{\mathbb R^\infty}t}u_0)$$
holds.

{\bf Proof:} 
About (2), we show it. About (1), it can be similarly shown. 
Let $f(t):=e^{\sqrt{-1}\triangle_{\mathbb R^N}t}f_0$ and 
$u(t):=e^{\sqrt{-1}\triangle_{\mathbb R^\infty}t}u_0$. 

First, we show that if $f_0\in H^2(\mathbb R^N)$ holds and 
$u_0$ is an element of the domain of $\triangle_{\mathbb R^\infty}$, 
then $e^{\sqrt{-1}\triangle_{\mathbb R^\infty}t}(f_0\otimes u_0)
=f(t)\otimes u(t)$ holds. From Lemma 43 and Lemma 38 (2), 
$$f\otimes u\in C^1((-\infty,+\infty);L^2(\mathbb R^\infty)),$$  
$$\frac{d}{dt}(f\otimes u)=\sqrt{-1}((\triangle_{\mathbb R^N}f)\otimes u
+f\otimes(\triangle_{\mathbb R^\infty}u))$$
hold. On the other hand, because from Lemma 38 (1), $f\otimes u\in H^1(\mathbb R^\infty)$ holds, 
from Lemma 42 (3) and Lemma 42 (4), for any $\varphi\in H^1(\mathbb R^\infty)$, 
$$\langle f\otimes u, \varphi\rangle_{L^2(\mathbb R^\infty)}-\langle (\triangle_{\mathbb R^N}f)\otimes u
+f\otimes(\triangle_{\mathbb R^\infty}u), \varphi\rangle_{L^2(\mathbb R^\infty)}
=\langle f\otimes u, \varphi\rangle_{H^1(\mathbb R^\infty)}
$$
holds. 
So, $f\otimes u-((\triangle_{\mathbb R^N}f)\otimes u
+f\otimes(\triangle_{\mathbb R^\infty}u))
=(1-\triangle_{\mathbb R^\infty})(f\otimes u)$ holds. 
$$(\triangle_{\mathbb R^N}f)\otimes u
+f\otimes(\triangle_{\mathbb R^\infty}u)
=\triangle_{\mathbb R^\infty}(f\otimes u)$$
holds. Therefore, if $f_0\in H^2(\mathbb R^N)$ holds and 
$u_0$ is an element of the domain of $\triangle_{\mathbb R^\infty}$, 
then $e^{\sqrt{-1}\triangle_{\mathbb R^\infty}t}(f_0\otimes u_0)
=f(t)\otimes u(t)$ holds. 

There exist a sequence $\{g_{0,n}\}_n$ in $H^2(\mathbb R^N)$ 
such that $\lim_n\|g_{0,n}-f_0\|_{L^2(\mathbb R^N)}=0$ holds. 
There exist a sequence $\{v_{0,n}\}_n$ in the domain of $\triangle_{\mathbb R^\infty}$ 
such that $\lim_n\|v_{0,n}-u_0\|_{L^2(\mathbb R^\infty)}=0$ holds. 
Let $g_n(t):=e^{\sqrt{-1}\triangle_{\mathbb R^N}t}g_{0,n}$ and 
$v_n(t):=e^{\sqrt{-1}\triangle_{\mathbb R^\infty}t}v_{0,n}$. 
Then, $\lim_n\|g_{n}(t)-f(t)\|_{L^2(\mathbb R^N)}=0$ 
and $\lim_n\|v_{n}(t)-u(t)\|_{L^2(\mathbb R^\infty)}=0$ hold. 
From Lemma 32, Lemma 33 and Lemma 38 (2), 
$$\lim_n\|g_n(t)\otimes v_{n}(t)-f(t)\otimes u(t)\|_{L^2(\mathbb R^\infty)}=0$$ 
holds. On the other hand, from 
$\lim_n\|g_{0,n}\otimes v_{0,n}-f_0\otimes u_0\|_{L^2(\mathbb R^\infty)}=0$, 
$$\lim_n\|e^{\sqrt{-1}\triangle_{\mathbb R^\infty}t}(g_{0,n}\otimes v_{0,n})
-e^{\sqrt{-1}\triangle_{\mathbb R^\infty}t}(f_0\otimes u_0)\|_{L^2(\mathbb R^\infty)}=0$$ 
holds. So, because of 
$e^{\sqrt{-1}\triangle_{\mathbb R^\infty}t}(g_{0,n}\otimes v_{0,n})
=g_n(t)\otimes v_{n}(t)$, 
it follows. 
\hfill 
$\blacksquare$

\newpage

\sect{Inseparability of $L^2(\mathbb R^\infty)$} 
{\bf Lemma 45:} 
Let $\{f_n\}_{n\in\mathbb N}$ be a sequence of $[0,+\infty)$-valued measurable functions on $\mathbb R$. 
Suppose that for any $n\in\mathbb N$, $\|f_n\|_{L^2(\mathbb R)}=1$ holds. 
Then, the followings hold. 

(1) Let $m\in \mathbb N$. Then, 
$$(\prod_{n\in\{m\}}f_n(x_n)^2dx_n)(\prod_{n\in\mathbb N\setminus\{m\}}f_n(x_n)^2dx_n)
=\prod_{n\in\mathbb N}f_n(x_n)^2dx_n$$
holds. 

(2) Let $\{a_n\}_{n\in\mathbb N}\in \mathbb R^\infty$. Then, 
$$\tau_{\{a_n\}_{n\in\mathbb N}}(\prod_{n\in\mathbb N}f_n(x_n)^2dx_n)
=\prod_{n\in\mathbb N}f_n(x_n-a_n)^2dx_n$$
holds.

{\bf Proof:} 
Although it is a natural result, we include the proof just in case. 

Let $N\in \mathbb N$. Let $\{E_n\}_{n=1}^N$ be a family 
of measurable sets of $\mathbb R$. 

(1)  
$((\prod_{n\in\{m\}}f_n(x_n)^2dx_n)(\prod_{n\in\mathbb N\setminus\{m\}}f_n(x_n)^2dx_n))
((\prod_{n=1}^NE_n)\times (\prod_{n=N+1}^\infty\mathbb R))
=\prod_{n=1}^N(\int_{E_n}f_n(x_n)^2dx_n)$ holds. 

(2) For $n\in \{1,2,\cdots,N\}$, let $F_{n}:=\{x_n\in\mathbb R|x_n+a_n\in E_n\}$. Then, 
$(\tau_{\{a_n\}_{n\in\mathbb N}}(\prod_{n\in\mathbb N}f_n(x_n)^2dx_n))
((\prod_{n=1}^NE_n)\times (\prod_{n=N+1}^\infty\mathbb R))
=(\prod_{n\in\mathbb N}f_n(x_n)^2dx_n)$


\noindent 
$((\prod_{n=1}^NF_n)\times (\prod_{n=N+1}^\infty\mathbb R))
=\prod_{n=1}^N(\int_{F_n}f_n(x_n)^2dx_n)
=\prod_{n=1}^N(\int_{E_n}f_n(x_n-a_n)^2dx_n)$ holds. 
\hfill 
$\blacksquare$

{\bf Lemma 46:} 
Let $\{f_n\}_{n\in\mathbb N}$ and $\{g_n\}_{n\in\mathbb N}$ 
be sequences of $[0,+\infty)$-valued measurable functions on $\mathbb R$. 
Suppose that for any $n\in\mathbb N$, $\|f_n\|_{L^2(\mathbb R)}=\|g_n\|_{L^2(\mathbb R)}=1$ holds. 
Suppose that there exists $m\in \mathbb N$ such that $\langle f_m, g_m\rangle_{L^2(\mathbb R)}=0$ holds. 
Then, $\langle \prod_{n\in\mathbb N}f_n(x_n)^2dx_n, 
\prod_{n\in\mathbb N}g_n(x_n)^2dx_n\rangle_{CM(\mathbb R^\infty)}=0$ holds. 

{\bf Proof:} 
From Lemma 45 (1) and Proposition 28 (1), it follows.
\hfill 
$\blacksquare$

{\bf Lemma 47:} 
Let $\{f_n\}_{n\in\mathbb N}$ be a sequence of $[0,+\infty)$-valued measurable functions on $\mathbb R$. 
Suppose that for any $n\in\mathbb N$, $\|f_n\|_{L^2(\mathbb R)}=1$ holds. 
Let $m\in \mathbb N$. Suppose that $f_m\in H^1(\mathbb R)$ holds. Then, 
$\prod_{n\in\mathbb N}f_n(x_n)^2dx_n\in H^1_m(\mathbb R^\infty)$ and 
$\frac{\partial}{\partial x_m}(\prod_{n\in\mathbb N}f_n(x_n)^2dx_n)
=(\prod_{n\in\{m\}}f^\prime_n(x_n)|f^\prime_n(x_n)|dx_n)\cdot
(\prod_{n\in\mathbb N\setminus\{m\}}f_n(x_n)^2dx_n)$ hold.

{\bf Proof:} 
There uniquely exists $e_m\in \mathbb R^\infty$ such that  
$e_{m,m}=1$ holds and for any $n\in\mathbb N\setminus\{m\}$, $e_{m,n}=0$ holds. 
From Lemma 45, for $h>0$, 
$\tau_{he_m}(\prod_{n\in\mathbb N}f_n(x_n)^2dx_n)
=(\prod_{n\in\{m\}}f_n(x_n-h)^2dx_n)(\prod_{n\in\mathbb N\setminus\{m\}}f_n(x_n)^2dx_n)$
holds. So, because from Proposition 28 (2), 
$\frac{\tau_{he_m}(\prod_{n\in\mathbb N}f_n(x_n)^2dx_n)-\prod_{n\in\mathbb N}f_n(x_n)^2dx_n}{h}
+(\prod_{n\in\{m\}}f^\prime_n(x_n)|f^\prime_n(x_n)|dx_n)\cdot
(\prod_{n\in\mathbb N\setminus\{m\}}f_n(x_n)^2dx_n)
=(\prod_{n\in\{m\}}(\frac{f_n(x_n-h)-f_n(x_n)}{h}+f^\prime_n(x_n))
|\frac{f_n(x_n-h)-f_n(x_n)}{h}+f^\prime_n(x_n)|dx_n)\cdot(\prod_{n\in\mathbb N\setminus\{m\}}f_n(x_n)^2dx_n)$ 
holds, from Proposition 28 (1), it follows. 
\hfill 
$\blacksquare$

{\bf Proposition 48:} 
There exists an orthonormal system $\{u_\tau\}_{\tau\in \prod_{n\in\mathbb N}\{0,1\}}$ 
of $L^2(\mathbb R^\infty)$. 

{\bf Proof:} 
There exists $f\in C^\infty(\mathbb R)$ such that 
$\int_{\mathbb R}|f|^2dx=1$, for any $x\in \mathbb R$, $f(x)\geq 0$ holds 
and for any $x\in \mathbb R\setminus (0,1)$, $f(x)=0$ holds. 
There exists a sequence $\{L_n\}_{n\in\mathbb N}$ in $(0,+\infty)$ 
such that 
$$\sum_{n\in\mathbb N}\frac{1}{L_n^4}<+\infty$$ 
holds. For $\tau\in \prod_{n\in\mathbb N}\{0,1\}$, let 
$$u_\tau:=\prod_{n\in\mathbb N}(\frac{1}{L_n}f(\frac{x_n}{L_n^2}-\tau(n)))^2dx_n.$$ 
Then, from Lemma 46, $\{u_\tau\}_{\tau\in \prod_{n\in\mathbb N}\{0,1\}}$ 
is an orthonormal system of $CM(\mathbb R^\infty)$. 
On the other hand, because from Lemma 47 and Proposition 28 (1), 
for any $\tau\in \prod_{n\in\mathbb N}\{0,1\}$ 
and any $n\in \mathbb N$, 
$$\|\frac{\partial u_\tau}{\partial x_n}\|_{CM(\mathbb R^\infty)}^2
=\frac{1}{L_n^4}\|f^\prime\|_{L^2(\mathbb R)}^2$$
holds, for any $\tau\in \prod_{n\in\mathbb N}\{0,1\}$, 
$u_\tau\in H^1(\mathbb R^\infty)$ holds. So, 
$\{u_\tau\}_{\tau\in \prod_{n\in\mathbb N}\{0,1\}}$ 
is an orthonormal system of $L^2(\mathbb R^\infty)$. 
\hfill 
$\blacksquare$

\newpage 

\sect{Translation invariance of $\triangle_{\mathbb R^\infty}$}
{\bf Lemma 49:} 
Let $a\in \mathbb R^\infty$ and $k\in \mathbb N$. 
Then, the followings hold. 

(1) Let $u\in H^1_k(\mathbb R^\infty)$. Then, $\tau_a u\in H^1_k(\mathbb R^\infty)$ 
and $\frac{\partial (\tau_a u)}{\partial x_k}=\tau_a(\frac{\partial u}{\partial x_k})$ hold. 

(2) Let $u_1, u_2\in H^1_k(\mathbb R^\infty)$. Then, 
$\langle \frac{\partial (\tau_a u_1)}{\partial x_k}, \frac{\partial (\tau_a u_2)}{\partial x_k}\rangle_{CM(\mathbb R^\infty)}
=\langle \frac{\partial u_1}{\partial x_k}, \frac{\partial u_2}{\partial x_k}\rangle_{CM(\mathbb R^\infty)}$
holds.

{\bf Proof:} 
(1) From $\tau_b\tau_c=\tau_{b+c}$ and Proposition 13, it follows. 

(2) From (1) and Proposition 13, it follows. 
\hfill 
$\blacksquare$

{\bf Lemma 50:} 
Let $a\in \mathbb R^\infty$. 
Then, the followings hold. 

(1) Let $u\in H^1(\mathbb R^\infty)$. Then, $\tau_a u\in H^1(\mathbb R^\infty)$ holds. 

(2) Let $u_1, u_2\in H^1(\mathbb R^\infty)$. Then, 
$\langle \tau_a u_1, \tau_a u_2\rangle_{H^1(\mathbb R^\infty)}
=\langle u_1, u_2\rangle_{H^1(\mathbb R^\infty)}$
holds. 

{\bf Proof:} 
(1) From Lemma 49 (1) and Proposition 13, it follows. 

(2) From Lemma 49 (2) and Proposition 13, it follows. 
\hfill 
$\blacksquare$

{\bf Lemma 51:} 
Let $a\in \mathbb R^\infty$ and $f_1, f_2\in L^2(\mathbb R^\infty)$. 
Then,  $\tau_a f_1, \tau_a f_2 \in L^2(\mathbb R^\infty)$ 
and $\langle \tau_a f_1, \tau_a f_2\rangle_{L^2(\mathbb R^\infty)}
=\langle f_1, f_2\rangle_{L^2(\mathbb R^\infty)}$ hold. 

{\bf Proof:} 
From Lemma 50 (1) and Proposition 13, it follows. 
\hfill 
$\blacksquare$

{\bf Proposition 52:} 
Let $a\in \mathbb R^\infty$. Let $u$ be an element of the domain of $\triangle_{\mathbb R^\infty}$. 
Then, $\tau_a u$ is an element of the domain of $\triangle_{\mathbb R^\infty}$ and 
$\triangle_{\mathbb R^\infty}(\tau_a u)=\tau_a(\triangle_{\mathbb R^\infty}u)$ holds. 

{\bf Proof:} 
From Lemma 50, Lemma 51 and Proposition 13, for any $\varphi\in H^1(\mathbb R^\infty)$, 
$\langle \tau_a u-\tau_a(\triangle_{\mathbb R^\infty}u),  \varphi\rangle_{L^2(\mathbb R^\infty)}
=\langle  (1-\triangle_{\mathbb R^\infty})u,  \tau_{-a}\varphi\rangle_{L^2(\mathbb R^\infty)}
=\langle  u,  \tau_{-a}\varphi\rangle_{H^1(\mathbb R^\infty)}$ 


\noindent 
$=\langle  \tau_a u,  \varphi\rangle_{H^1(\mathbb R^\infty)}$ 
holds. So, $\tau_a u=(1-\triangle_{\mathbb R^\infty})^{-1}
(\tau_a u-\tau_a(\triangle_{\mathbb R^\infty}u))$ holds. 
\hfill 
$\blacksquare$

{\bf Remark:} 
(1) Gross ([4]) considered an analog on infinite-dimensional real Hilbert space of the Laplacian.  
Of course, $\mathbb R^\infty$ is not Hilbert space. 
It is translation invariant. 
However, it is not a self-adjoint operator.

(2) For $p\in[1,+\infty)$, let $L_p$ be Ornstein-Uhlenbeck operator in $L^p(\mu)$ with respect to Wiener measure $\mu$. 
Of course, it may be more appropriate to think of $L_p$ as an analog of Ornstein-Uhlenbeck operator on $\mathbb R^n$, 
rather than as an analog of Laplace operator on $\mathbb R^n$. 
$L_p$ is a generator of a contraction semigroup and one of the main characters in Malliavin calculus (e.g., [10]). 
$L_2$ is a self-adjoint operator corresponding to infinite-dimensional Dirichlet form (e.g., [9]). 
However, $L_p$ is not translation invariant. 
\hfill 
$\square$

\newpage

{\bf Smooth domain:} 
Let $\rho\in C^\infty_0(\mathbb R)$. Suppose that for any $t\in \mathbb R$, $0<\sum_{k=2}^\infty |\rho(t-(-1)^kk)|$ holds. Then, 
$$\{ \, \{x_k\}_k \in \mathbb R^\infty \, | \, 0<\sum_{k=2}^\infty \rho(x_1-(-1)^kk)x_k \, \}$$ 
is a non-trivial smooth domain of $\mathbb R^\infty$.  
\hfill 
$\square$

{\bf Gradient flow:} 
Let $\Omega$ be a compact Riemannian manifold. 
Let $M$ be a submanifold of Hilbert space $L^2(d\omega_\Omega)$. Then, 
the real inner product $\mathrm{Re}(\langle \cdot, \cdot \rangle_{L^2(d\omega_\Omega)})$ 
induces Riemannian metric on $M$.
So, we can think of the gradient flow on $M$ of the function 
$u \mapsto \frac{1}{2}\int_{\Omega} |\mathrm{grad}_\Omega u|^2d\omega_\Omega$. 
When $M=L^2(d\omega_\Omega)$ holds, it is the heat equation $u_t=\triangle_\Omega u$. 
\hfill 
$\square$

{\bf Impression:} 
For the measurable space $\mathbb R^n$, the set of all vectors corresponding to quantum states 
that are subject to probability interpretation is $CM(\mathbb R^n)\setminus\{0\}$, 
but on the other hand, the set of all vectors that can be solved by the ordinary Schr\"{o}dinger equation 
is $L^2(\mathbb R^n)$. It may be an ideal model that is easy to understand, 
rather than a realistic model that is difficult to understand. In particular, the difference between thinking 
that there exists the quantum state corresponding to the Dirac measure or not is essential. 
\hfill 
$\square$

\vfill









\newpage

{\it References} 

[1] 
Ay, N., Jost, J., L\^{e}, H. and Schwachh\"{o}fer, L., 
{\it Information geometry}, Springer, Cham, 2017. 

[2]
Brezis, H, {\it Functional analysis, Sobolev spaces and partial differential equations}, 
Springer, New York, 2011.

[3] 
Dunford, N. and Schwartz, J. T., 
{\it Linear operators}, John Wiley \& Sons, New York, 1988. 

[4] 
Gross, L., 
{\it Potential theory on Hilbert space}, 
J. Functional Analysis, 
1 (1967), 123--181. 

[5]
Halmos, P. R., {\it Measure Theory}, D. Van Nostrand Co., New York, 1950.

[6] 
H\"{o}rmander, L., 
{\it Fourier integral operators. I}, 
Acta Math., 127 (1971), 79--183. 

[7] 
Jost, J., {\it Postmodern analysis}, Springer-Verlag, Berlin, 2005.

[8] 
Kato, T, {\it Perturbation theory for linear operators}, Springer-Verlag, Berlin, 1995

[9] 
Kusuoka, S., 
{\it Dirichlet forms and diffusion processes on Banach spaces}, 
J. Fac. Sci. Univ. Tokyo Sect. IA Math., 
29 (1982), 79--95. 

[10] 
Malliavin, P., 
{\it Stochastic calculus of variation and hypoelliptic operators}, 
Proceedings of the International Symposium on Stochastic Differential Equations, 
pp. 195--263, Wiley, New York-Chichester-Brisbane, 1978. 

[11] 
Neveu, J., {\it Bases math\'{e}matiques du calcul des probabilit\'{e}s}, 
Deuxi\`{e}me \'{e}dition, Masson et Cie, Paris, 1970. 

[12]
Reed, M. and Simon, B., {\it Methods of Modern Mathematical Physics}, 
Academic Press, New York, 1975, 1979, 1978, 1980. 

\end{document}